\documentclass[12pt]{amsart}
\usepackage{amsmath}
\usepackage{amssymb}
\usepackage{mathrsfs}       
\usepackage[curve]{xypic}

\usepackage[colorlinks=true, urlcolor=blue,bookmarks=true,bookmarksopen=true, citecolor=blue,hypertex]{hyperref}

\numberwithin{equation}{section}

\newtheorem{theorem}{Theorem}[subsection]
\newtheorem{defn}[theorem]{Definition}
\newtheorem{assumption}[theorem]{Assumption}

\newtheorem{corollary}[theorem]{Corollary}
\newtheorem{example}[theorem]{Example}

\newtheorem{lemma}[theorem]{Lemma}
\newtheorem{prop}[theorem]{Proposition}
\newtheorem{remark}[theorem]{Remark}

\def \begineq{\begin{equation}}
\def \endeq{\end{equation}}

\def \bb{\mathbb}

\def \mc{\mathcal}
\def \mf{\mathfrak}
\def \ms{\mathscr}

\def \CC{{\bb{C}}}

\def \JJ{{\bb{J}}}

\def \NN{{\bb{N}}}

\def \RR{{\bb{R}}}
\def \TT{{\bb{T}}}
\def \UU{{\bb{U}}}
\def \VV{{\bb{V}}}
\def \WW{{\bb{W}}}
\def \ZZ{{\bb{Z}}}

\def \DDC{{\mc D}}

\def \JJC{{\mc{J}}}
\def \KKC{{\mc{K}}}
\def \LLC{{\mc{L}}}

\def \TTC{{\mc T}}

\def \GGS{{\ms G}}

\def \({\left(}
\def \){\right)}
\def \<{\langle}
\def \>{\rangle}

\def \dsum{\oplus}

\def \inter{\cap}

\def \tensor{\otimes}
\def \union{\cup}
\def \vargeq{\geqslant}

\def \xto{\xrightarrow}

\def \Ann{{\rm Ann}}
\def \Aut{{\rm Aut}}

\def \Diff{{\rm Diff}}

\def \Hom{{\rm Hom}}

\def \Span{{\rm Span}}

\def \qed{\hfill $\square$ \vspace{0.03in}}

\renewcommand{\1}{1\!\!1}
\def \ttriangle{{\triangle \!\!\!\!\!\triangle}}

\begin{document}

\title{Extended manifolds and extended equivariant cohomology}

\author{Shengda Hu}

\address{D\'epartement de Math\'ematiques et de Statistique, Universit\'e de  Montr\'eal, CP 6128 succ Centre-Ville, Montr\'eal, QC H3C 3J7, Canada}
\email{shengda@dms.umontreal.ca}
\author{Bernardo Uribe}

\address{Department of mathematics, University of Los Andes, Carrera 1 n. 18A -10, Bogot\'a, Colombia}
\email{buribe@uniandes.edu.co}
\thanks{The first author was supported by funding from DMS at Universit\'e de Montr\'eal. Part of the research of the first author was carried out at the Yantze Center of Mathematics at Sichuan University. The second author was partially supported by the  ``Fondo de apoyo
a investigadores jovenes" from Universidad de los Andes.
 Part of the research of the second author was carried out at the Max Planck Institut in Bonn
 the Erwin Schr\"odinger Institut in Vienna.  }



\abstract We define the category of manifolds with extended
tangent bundles, we study their symmetries and we consider the
analogue of equivariant cohomology for actions of Lie groups in
this category. We show that when the action preserves the
splitting of the extended tangent bundle, our definition of
extended equivariant cohomology agrees with the twisted
equivariant de Rham model of Cartan, and for this case we show
that there is localization at the fixed point set, \`a la
Atiyah-Bott.
\endabstract
\maketitle

\section{Introduction}\label{intro}
The study of the geometry on the generalized tangent bundle $\TT M
= TM \dsum T^*M$ starts with the paper \cite{Hitchin}, in which
Hitchin introduced the notion of generalized complex structures.
The framework of generalized complex geometry was first developed
by Gualtieri in his thesis \cite{Gualtieri} and much more has been
done since.
The algebraic structure underlying the considerations of
generalized geometry is the structure of a Courant algebroid on
$\TT M$ (cf. definition \ref{mfdext:courant}), this is the main
object of study in this article.

The natural  group of symmetries  for the tangent bundle $TM$ is
the diffeomorphism group $\Diff(M)$. In contrast, the natural
group of symmetries for $\TT M$ is given by $\ms G = \Diff(M)
\ltimes \Omega^2(M)$. The action of $B \in \Omega^2(M)$, which is
also called a $B$-transformation, is given by
$$e^B (X + \xi) = X + \xi + \iota_X B, \text{ for } X \in \Gamma(TM) \text{ and }
 \xi \in \Omega^1(M).$$
In this article, we take this approach one step further and take
the following point of view. We see the $B$-transformations as
part of an \emph{extended change of coordinates}. As in classical
differential geometry that intrinsic quantities do not depend on
the choice of coordinates, here we ask for independence on the
choice of extended coordinates.

With this point of view, the splitting of $\TT M$ into direct sum
is only a choice of extended coordinates, in which applying a
$B$-transformation corresponds to a different choice of
coordinates. Along the same line, the embedding of $T^* M \to \TT
M$ is invariant with respect to the extended change of
coordinates. Thus we are led to the consideration of the extension
sequence:
$$0 \to T^*M \to \TTC M \to TM \to 0,$$
where $\TTC M$ is the same bundle $\TT M$ but without a preferred
splitting. The Courant algebroid $\TTC M$ of this type is called
\emph{exact} \cite{Severa} and we will use the name \emph{extended
tangent bundle} to emphasize the relation with the geometry of the
manifold, as well as the absence of a preferred splitting. In the
following, we will use the notion \emph{extended manifold} to
denote a manifold $M$ with an extended tangent bundle $\TTC M$.

This paper accomplishes the following. In the first part
(\S\ref{linearext} and \ref{mfdext}), starting with linear algebra
(\S\ref{linearext}) we describe the category of extended manifolds
(\S\ref{mfdext}). In the second part (\S\ref{symmetry} and
\ref{equicoho}), we consider the induced action on the spinors
(\S\ref{symmetry}) and describe a new equivariant cohomology
(\S\ref{equicoho}) which involves non-trivially the
\emph{extended} part of the symmetries, i.e. the action of
$\Omega^2(M)$. In \S\ref{equicoho} we also consider some special
cases and examples. In the last part (\S\ref{kirwan} and
\ref{local}), we check the properties of the equivariant
cohomology we introduced (\S\ref{kirwan}) and show that, in
certain cases, the new equivariant cohomology allows localization
to fixed point sets (\S\ref{local}), in the sense of Atiyah-Bott.

In the first part, to define the category, we are mainly concerned
with the definition of morphisms (cf. definition
\ref{mfdext:morphdef}). There are two ingredients in the definition,
\emph{reversed} structure (cf. definition \ref{mfdext:reversion}) and
\emph{isotropic extended submanifold} (cf. definition
\ref{mfdext:extendedseq}).
Recall that the structure of a Courant algebroid on $\TTC M$ is
given by the datum $(*, \<,\>, a)$ ( cf. definition
\ref{mfdext:courant}), where $\<,\>$ is a non-degenerate symmetric
pairing. The \emph{reversed} structure of Courant algebroid on
$\TTC M$ is then given by the datum $(*, -\<,\>, a)$, which we
denote by $-\TTC M$. A standard reversion is by definition an
isomorphism of Courant algebroids $\tau : \TTC M \to - \TTC M$,
covering the identity on $TM$. Standard reversions can be seen as
a more intrinsic way of representing isotropic splittings of $\TTC
M$, as they are in one-to-one correspondence to each other. Then a
morphism between $(M, \TTC M)$ and $(N, \TTC N)$ is a map $f : M
\to N$ whose graph is an isotropic extended submanifold of $(M
\times N, -\TTC M \boxplus \TTC N)$. We show that this naturally
encodes the action of $\Omega^2(M)$. Summarizing, we have:

\vspace{0.1in}
 \emph{
The category $\mc ESmth$ of extended manifolds is given by the
following. An object in the category is a manifold $M$ together
with an extended tangent bundle $\TTC M$. An morphism $\tilde f :
(M, \TTC M) \to (N, \TTC N)$ is given by a smooth map $f: M \to N$
and an isotropic extended structure $\mc E$ on the graph of $f$,
with respect to the extended tangent bundle $-\TTC M \boxplus \TTC
N$. \\ When we choose and fix splittings of the extended tangent
bundles, a morphism is given equivalently by a smooth map $f : M
\to N$ and a two form $b_{\tilde f} \in \Omega^2(M)$ so that $H_M
= db_{\tilde f} + f^*H_N$, where $H_\cdot$ are the twisting forms
defined by the chosen splittings and $[H_\cdot]$ give the \v
Severa classes of the extended tangent bundles. In this way, the
composition of two morphisms $\tilde f = (f, b_{\tilde f})$ and
$\tilde g = (g, b_{\tilde g})$ is given by $\tilde h = (g\circ f,
b_{\tilde f} + f^*b_{\tilde g})$. } \vspace{0.1in}

\noindent Although the definition of composition as above uses the
splittings of extended tangent bundles, we show that the morphism
thus obtained does not depend on the choice of splitting.

The second part discusses equivariant cohomology, in the extended
setting. Recall that the space of spinors for $\TTC M$ with
pairing $\<,\>$ can be identified with $\Omega^\bullet(M)$. From
our point of view, the identification is not canonical, rather, it
depends on the choice of extended coordinates. We thus arrive at
the notion of \emph{abstract spinor space} $\mc S^\bullet(M)$ (cf.
definition \ref{spinor:abstract}). For $\mf g$ be a Lie algebra,
we introduce two notions of $\mf g$-actions on $\TTC M$. The first
is the \emph{generalized action} of $\mf g$, given by a Lie
algebra homomorphism $\mf g \to \ms X_\TTC$, where $\ms X_\TTC$ is
the Lie algebra of infinitesimal symmetries of $\TTC M$. The
second is when the generalized action factors through the
composition $\mf g \stackrel{\delta}{\to} \Gamma(\TTC M)
\stackrel{\kappa}{\to} \ms X_\TTC$, where the map $\kappa$ sends
the bracket of the Courant algebroid $\TTC M$ to the Lie bracket
on $\ms X_\TTC$. We say then that this generalized action is an
\emph{extended action}  and
we have obtained the following result.

\vspace{0.1in}
\emph{
Choose and fix a splitting $s$ of $\TTC M$.
For $\mf X = s(X) + \xi \in \Gamma(\TTC M)$, let $\LLC^\TTC_{\mf X}
= \LLC_X + (d\xi - \iota_X H) \wedge$, $\iota_{\mf X} = \iota_X + \xi \wedge$
and $d_\TTC = d - H \wedge$, where $H$ is the twisting form defined by $s$.
 Then these operators are independent of the choice of splitting $s$.
  Let $\delta : \mf g \to \Gamma(\TTC M)$ be an extended $\mf g$-action that is isotropic.
   Then the following is a chain complex
$$C^\bullet_{\mf g}(\TTC M) = \left\{\rho \in \mc S^\bullet(M) \tensor \widehat{S}(\mf g^*)
\left|\LLC^\TTC_{\delta(\tau)} \rho = 0  \ \mbox{for all} \ \tau
\in \mf g \right.\right\}, \,\, d_{\TTC, \delta} = d_\TTC - \sum_j
u_j \iota_{\delta(\tau_j)}.$$  The cohomology $H^\bullet_{G}(\TTC
M)$ of $(C^\bullet_{\mf g}(\TTC M), d_{\TTC , \delta })$ is called
the \emph{extended $\mf g$-equivariant cohomology}. }
\vspace{0.1in}

\noindent We note that even when the generalized action of $\mf g$
on $\TTC M$ is trivial, the map $\delta : \mf g \to \Gamma(\TTC
M)$ may not be trivial and the above cohomology is different from
the tensor product of $H^*(M)$ and $\widehat{S}(\mf g^*)$. In fact, as the
examples in \S\ref{equicoho}
show, the extended equivariant cohomology can be genuinely
different from the ordinary equivariant cohomology.

The last part concerns the case when the $\mf g$-action is
integrable to an action of compact Lie group $G$. In this case,
let $F$ be the fixed point loci. Then the induced extended tangent
bundle $\TTC F$ naturally embeds into $\TTC M$. The naturality of
this embedding does not hold for general submanifolds. When the
$\mf g$-action is an extended action that preserves a splitting,
we show that the extended $\mf g$-equivariant cohomology localizes
to the fixed point set in the sense of Atiyah-Bott. Along the way,
we show that the axioms of cohomology theory are satisfied by the
extended equivariant cohomology and there is a Thom isomorphism
between the cohomology of the base and the vertically compactly
supported cohomology of a vector bundle. The later is the corner
stone for the localization argument. \vspace{0.1in}

We note that the terminology of \emph{extended action} has been
introduced by \cite{Bursztyn} in a slightly different fashion, in
connection with their reduction construction.
Various authors have described the twisted cohomology
$H^\bullet(M; H)$ in detail \cite{Mathai, Atiyah-Segal} and the
idea of twisting a equivariant cohomology theory by a $3$-class is
also known \cite{FreedHopkinsTeleman}. As far as the authors know,
the extended complex introduced in this article for the extended
actions has not been discussed before in the literature.

\vspace{0.1in} \noindent {\bf Acknowledgements.} The authors would
like to thank H. Bursztyn, A. Cardona, G. Cavalcanti, R. Cohen, M. Crainic, N.
Hitchin, K. Hori, E. Lupercio, V. Mathai, C. Teleman, for many
fruitful conversations. Also the authors would like to thank the
support of the Mathematical Sciences Research Institute, the Max
Planck Institut  and the Erwin Schr\"odinger Institut.


\section{Extended linear algebra}\label{linearext}
\subsection{Extended linear spaces}\label{linearext:space}
Let $V$ be a ($\RR$-)linear space of dimension $n$. An
\emph{extension} of $V$ is a triple $(\VV, \hat V, \<,\>)$ where
$\VV$ is a ($\RR$-)linear space of dimension $2n$, with a
non-degenerate pairing $\<,\>$ of signature $(n,n)$, so that $\hat
V \subset \VV$ is a maximal isotropic subspace which fits into the
extension sequence
\begin{equation}\label{linearext:extseq} 0 \to \hat V \xto{a^*} \VV \xto{a} V \to 0, \end{equation}
where $\hat V \simeq V^*$ via the pairing $2\<,\>$.
For $\mf X \in \VV$, we write $X = a(\mf X)$, then we have $\xi(X) = 2\<\hat\xi, \mf X\>$ for any $\xi \in V^*$, where we use the ``hat''  $\hat{} : V^* \to \hat V : \xi \mapsto \hat \xi$ to denote the identification induced by the pairing.

A splitting $s : V \to \VV$ of the above exact sequence is called \emph{isotropic} if the image is isotropic with respect to the pairing.
The set $\mc I$ of isotropic splittings $s$ is a torsor over $\wedge^2 V^*$. In fact, let $B(X, Y) = \<(s'-s)(X), Y\>$, then
$$\<(s'-s)(X), Y\> = -\<X-s'(X), s'(Y)\> + \<X-s(X), s(Y)\> = -\<X, (s'-s)(Y)\>,$$
i.e. $s'(X) = (B\circ s)(X) := s(X) + \widehat{\iota_X B}$.

\subsection{Extended subspaces}\label{linearext:subspace}
Let $i: W \subset V$ be a linear subspace.
An \emph{extended structure} $E$ on $W$ is a subspace of $\VV$ which fits
into the following diagram:
\begin{equation}\label{linearext:extstr}0 \to C \to E \xto a W \to 0, \text{ for some } C \subset \hat V,
\end{equation}
where $a$ is the restriction of the projection $\VV \to V$.
The \emph{extended subspace} $\mc W = (W, E)$ is \emph{isotropic}
 (resp. \emph{non-degenerate}) if the restriction of $\<,\>$ to $E$
 vanishes (resp. is non-degenerate). In the following, we will
 only consider \emph{maximal} extended structures, i.e. those that are
not contained in any other extended structures of the same type.
Let $K = \Ann_{\hat V} W$, $\hat W = \hat V / K$ and $\WW =
\Ann(K)/ K$, then we have an induced extension of $W$ given by
$(\WW, \hat W, \<,\>_W)$.

\subsubsection{Isotropic}\label{linearext:isosub} Let $E$ be a (maximal --- we will drop maximal in the following, as they are the only kind we consider here) isotropic extended structure on $W$, then $C = K$.
The extension sequence \eqref{linearext:extstr} descends to $\WW$ and becomes $0 \to E/K \to W \to 0$, i.e. $E/K$ is an isotropic splitting of $\WW$. On the other hand, for any splitting $s_W$ of $\WW$, the preimage of $s_W(W)$ in $\Ann(K)$ under the quotient map gives an isotropic extended structure on $W$. It follows that the space of isotropic extended structures on $W$ is a torsor over $\wedge^2 W^*$.

%

\subsection{Reversion} \label{linearext:reversion} Let $(\VV, \hat V, \<,\>)$
be an extension of $V$, then the \emph{reversed} extension is
defined to be $(\VV, \hat V, -\<,\>)$. When we abbreviate $\VV$ as
the extension, the reversed extension will be denoted $-\VV$. A
\emph{standard reversion} of $\VV$ is an isomorphism $\tau: \VV
\to -\VV$ which covers the identity on $V$. Let $s : V \to \VV$ be
an isotropic splitting, then we define the standard reversion
$\tau_s : \VV \to -\VV : \mf X = s(X) + \hat \xi \mapsto s(X) -
\hat \xi$, since
$$\<\tau_s(\mf X), \tau_s(\mf Y)\> = \<s(X) - \hat \xi, s(Y)-\hat \eta\> = \frac{1}{2}(\<\xi, Y\> + \<\eta, X\>) = - \<\mf X, \mf Y\>.$$
On the other hand, for any standard reversion, the sub-space
$V_\tau$ fixed by $\tau$ is non-empty and maximally isotropic,
since $\<\mf X, \mf Y\> = 0$ for all $\mf X, \mf Y \in V_\tau$. In
fact, we see that splittings and standard reversions are in
one-to-one correspond to each other, i.e. the set of standard
reversions is a torsor over $\wedge^2 V^*$.

\subsection{Extended complex structures}\label{linearext:gencplx}
A \emph{linear extended complex structure}
 on $\VV$ is defined as a linear map $\JJ : \VV \to \VV$ so that $\JJ^2 = -\1$ and
$\<\JJ \cdot, \JJ \bullet\> = \<\cdot, \bullet\>$. Let $P : V^*
\simeq \hat V \to \VV \xto{\JJ} \VV \to V$ be the composition of
the natural maps, then from the definition of $\JJ$ we have $P^* =
-P$ . It follows that $P \in \wedge^2 V$. An \emph{extended
Lagrangian} (resp. \emph{complex}) subspace $W \subset V$ is an
isotropic (resp. non-degenerate) extended subspace $(W, E)$ so
that $E$ is preserved by $\JJ$. When $W$ is extended Lagrangian,
$K = \Ann_{\hat V} W$ is isotropic with respect to $P$, since we
have the restriction $P: K \to W = \Ann_V K$, correspondingly, we
say that $W$ is \emph{co-isotropic} with respect to $P$. We may
also see $\JJ$ as a linear extended complex structure on $-\VV$ as
well. Let $\tau : \VV \to -\VV$ be a standard reversion, then
$\JJ_\tau = \tau \circ \JJ \circ \tau^{-1} : \VV \to \VV$ is the
\emph{$\tau$-reversed} linear extended complex structure (which
corresponds to the \emph{twisted} structure in \cite{Ben-Bassat}).

\subsection{Direct sums and morphisms}\label{linearext:morph}
Let
$(\VV, \hat V, \<,\>_V)$ and $(\WW, \hat W, \<,\>_W)$ be two extended linear spaces. The \emph{direct sum} extension of $V \dsum W$ is defined as $(\VV \dsum \WW, \hat V\dsum \hat W, \<,\>_V \dsum \<,\>_W)$. Let $\phi : W \to V$ be a linear map. An \emph{extension} $\tilde \phi = (\phi, E)$ of $\phi$ is given by an isotropic extended structure $E \subset -\WW \dsum \VV$ over $graph(\phi) \subset W \dsum V$,
i.e.
$$0 \to K_\phi \to E \to graph(\phi) \to 0.$$
By the arguments in \S\ref{linearext:isosub} we see that $(\hat
\eta, \hat \xi) \in K_\phi \subset \hat W \dsum \hat V$ is given
by
$$-\<\hat \eta, w\>_W + \<\hat \xi, \phi(w)\>_V = -\eta(w) + \xi(\phi(w)) = 0 \text{ for all } w \in W \Rightarrow K_\phi = graph(\hat{\phi}),$$
where $\phi^* : V^* \to W^*$ is dual to $\phi$, so that
$\hat{\phi} : \hat V \to \hat W$ via the identification
$\hat{\left.\right.}$.
\begin{defn}\label{linearext:morphdef}
We call $\tilde \phi$ a \emph{morphism} between the extended linear
 spaces and denote it by $\tilde \phi = (\phi, E) : \WW \to \VV$.
\end{defn}

Let's look at the subspace $E$ in more detail. Choose splittings $s_V$ and $s_W$
 of $\VV$ and $\WW$ respectively. Then an element $e \in E$ is of the form
\begin{equation} \label{linearext:descr}
\begin{split}
e = & \hat \phi(\hat \eta) + \hat \eta + s_W(x) + s_V(\phi(x)) + \lambda_{\hat V} (x) + \lambda_{\hat W}(x)\\
= & \hat \phi(\hat \eta') + \hat \eta' + s_W(x) + s_V(\phi(x)) + (\lambda_{\hat W} - \hat \phi \circ \lambda_{\hat V})(x)\\
\end{split}
\end{equation}
$$\text{ where } x \in W, \,\, \hat \eta \text{ and } \hat\eta' = \hat \eta + \lambda_{\hat V}(x)\in \hat V, \,\, \lambda_{\hat W} : W \to \hat W \text{ and } \lambda_{\hat V} : W \to \hat V.$$
We regard a linear map $q : W \to \hat W$ also as a quadratic form via $q(x, x') = 2\<q(x), s_W(x')\>_W$, where $s_W$ is any isotropic splitting of $\WW$. We then check that $E$ is isotropic iff
$$b_E := \lambda_{\hat W} - \hat \phi \circ \lambda_{\hat V} \in \wedge^2 W^*.$$
It follows that $E$ determines, and is uniquely determined, by
$b_E \in \wedge^2 W^*$ and we obtain the explicit form of the
statement in \S\ref{linearext:isosub}.
\begin{prop}\label{linearext:morphform}
Let $\phi : W \to V$ be a linear map. Choose splittings $s_W$ and $s_V$ of $\WW$
 and $\VV$ respectively, then any morphism $\tilde \phi : \WW \to \VV$ extending
 $\phi$ is determined by an element $b_{\tilde\phi} \in \wedge^2 W^*$.
  Let $B_W \in \wedge^2 W^*$ and $B_V \in \wedge^2 V^*$ and consider
  splittings $s'_W = s_W + \iota_\bullet B_W$ and $s'_V = s_V + \iota_\bullet B_V$ . Then the element
  $b'_{\tilde \phi} = b_{\tilde \phi} - B_W + \phi^* B_V \in \wedge^2 W^*$
   determines the same morphism $\tilde \phi$.
\end{prop}
{\it Proof:}
Directly compare the spaces
$$E = \{s_W(x) + \hat\phi(\eta) + \widehat{\iota_X b_{\tilde \phi}} + s_V(\phi(X)) + \eta | x \in W \text{ and } \eta \in \hat V\} \text{ with }$$
$$E' = \{s'_W(x) + \hat\phi(\eta) + \widehat{\iota_X b'_{\tilde \phi}} + s'_V(\phi(X)) + \eta | x \in W \text{ and } \eta \in \hat V\}$$
\qed

We note that we may always choose splittings of $\VV$ and $\WW$ so
that the element $b_{\tilde \phi}$ in the proposition
\ref{linearext:morphform} vanishes. In this case, we see that
$\tilde \phi$ is given by the maps $\phi : W \to V$ and $\phi^* :
V^* \to W^*$ on $\WW = W \dsum W^*$ and $\VV = V\dsum V^*$, where
the identifications come from the splittings as well as $\hat {}$.

\subsection{Composition}\label{linearext:comp}
Let $\tilde \phi = (\phi, E) : \WW \to \VV$ and $\tilde \psi = (\psi, F) : \VV \to \UU$ be two morphisms of extended linear spaces. We consider the composition $\lambda = \psi \circ \phi : W \to U$ and the induced extended structure on $graph(\lambda)$. According to the description given in the previous subsection, we choose and fix splittings of $\WW$, $\VV$ and $\UU$, then the morphism can be represented by the pair $(\phi, b_{\tilde \phi})$ and $(\psi, b_{\tilde \psi})$, where $b_{\tilde \phi} \in \wedge^2 W^*$, for example. We then define the composition of the morphisms:
\begin{defn}\label{linearext:composit}
With the chosen splittings of the extended linear spaces, the \emph{composition} $\tilde \lambda = \tilde \psi \circ \tilde \phi$ is given by $\lambda = \psi \circ \phi$ and $b_{\tilde \lambda} = b_{\tilde \phi} + \phi^* b_{\tilde \psi} \in \wedge^2 W^*$.
\end{defn}
We check that the composition is well-defined, i.e. it does not depend
on the choice of splittings, for which we use the proposition \ref{linearext:morphform}.
 Let $B_W \in \wedge^2 W^*$, $B_V \in \wedge^2 V^*$ and $B_U \in \wedge^2 U^*$ defining different splittings of the extended
 linear spaces. Then the same morphisms $\tilde \phi$ and $\tilde \psi$ are given by
$$b'_{\tilde \phi} = b_{\tilde \phi} -  B_W + \phi^* B_V \text{ and } b'_{\tilde \psi} = b_{\tilde \psi} - B_V + \psi^* B_U,$$
and the definition then gives
$$b'_{\tilde \lambda} = b'_{\tilde \phi} + \phi^* b'_{\tilde \psi} = b_{\tilde \phi} + \phi^* b_{\tilde \psi} - B_W + \lambda^* B_U =b_{\tilde \lambda}- B_W + \lambda^* B_U,$$
which by proposition \ref{linearext:morphform} gives the same
morphism $\tilde \lambda$. It is  easy to check that the following
equality of compositions of extended morphisms holds
$$\tilde\phi_1 \circ (\tilde \phi_2 \circ \tilde \phi_3) = (\tilde \phi_1 \circ \tilde \phi_2) \circ \tilde \phi_3.$$

\subsection{The category}\label{linearext:cat}
Now the category $\mc EVect$ of extended linear spaces is defined. An object in $\mc EVect$ is given by the triple $(\VV, \hat V, \<,\>)$ in the exact sequence
$$0 \to \hat V \to \VV \to V \to 0,$$
so that $\<,\>$ is symmetric, non-degenerate of split signature and $\hat V$
 is maximally isotropic. We will usually denote such an object by $\VV$. A morphism $\tilde \phi : \WW \to \VV$ in $\mc EVect$ is given by an isotropic extended structure $E \subset -\WW \dsum \VV$ on the graph of $\phi \in \Hom(W, V)$. We will call $\tilde \phi$ to be an extension of $\phi$. The composition of morphisms is given by definition \ref{linearext:composit}.
In particular, from the definitions, we have the following fibration sequence:
$$\wedge^2 W^* \to \Hom_{\mc EVect}(\WW, \VV) \to \Hom(W, V).$$
When $\WW = \VV$, we may consider the group of invertible morphisms $\Aut_{\mc EVect}(\VV)$, which fits into the extension sequence
$$0 \to \wedge^2 V^* \to \Aut_{\mc EVect} (\VV) \to GL(V) \to 1.$$
The identity element in $\Aut_{\mc EVect}(\VV)$ is given by the diagonal $\ttriangle \subset -\VV \dsum \VV$, which is a morphism extending $id \in GL(V)$.


\subsection{Spinors}\label{linearext:spinor}
Let $Cl(\VV)$ be the Clifford algebra of $\VV$, i.e.
$$Cl(\VV) = \tensor^* \VV / (\mf X \otimes  \mf X - \<\mf X, \mf X\>).$$
For the choice of an isotropic splitting, we determine an action
of $Cl(\VV)$ on $\wedge^* \hat V$, exhibiting it as a spinor space
for $Cl(\VV)$. For $\rho \in \wedge^*\hat V$, we denote the action
of ${\mf X}$ on $\rho$ by the contraction
$$\iota_{\mf X} \rho = \iota_X \rho + \hat\xi \wedge \rho, \text{ where } \hat\xi = \mf X - s(X).$$

Let $B \in \wedge^2 V^*$, then $e^B = \sum_{j} \frac{B^{\wedge j}}{j!} \in \wedge^{ev} V^*$ and the action of $\wedge^2 V^*$ on $\wedge^* \hat V$ is defined by:
$$B \circ \rho = e^{-\hat B} \wedge \rho.$$
\begin{defn}\label{linearext:spinorabsp}
The \emph{abstract spinor space} $\mc S$ for the triple $(\VV, \hat V, \<,\>)$ is given by $\mc I \times_{\wedge^2 V^*} \wedge^* \hat V$, where the action is the anti-diagonal action (as usual).
\end{defn}
The bundle $\mc I \times \wedge^* \hat V$ over $\mc S$ with fiber $\wedge^2 V^*$
provides identification of $\mc S$ with the section $\{s\} \times \wedge^* \hat V$
 upon a choice of splitting $s \in \mc I$. The representation of $Cl(\VV)$
  on $\mc S$ is defined by any of such identifications as we can
  see in the following set of equations:
\begin{equation*}
\begin{split}
& e^{-\hat B} \wedge \left(\iota_X (e^{\hat B} \wedge \rho) + (\mf X - s(X) - \widehat {\iota_X B}) \wedge e^{\hat B} \wedge \rho\right)\\
= & e^{-\hat B}\wedge \left(\widehat{\iota_X B} \wedge e^{\hat B} \wedge \rho + e^{\hat B} \wedge \iota_X \rho + (\mf X - s(X) - \widehat{\iota_X B}) \wedge e^{\hat B} \wedge \rho\right)\\
= & \iota_X \rho + (\mf X - s(X)) \wedge \rho.
\end{split}
\end{equation*}
We will use $\iota_{\mf X}$ to denote the above action of $\mf X \in \VV$ on $\mc S$, which generates the $Cl(\VV)$-action.

\begin{lemma}\label{linearext:pullback}
Let $\tilde \phi : \WW \to \VV$ be a morphism of extended linear spaces. Then it induces natural \emph{pull-back} map of spinor spaces: $\tilde \phi^\bullet : \mc S_V \to \mc S_W$.
\end{lemma}
{\it Proof:}
We use the description of morphisms given in proposition \ref{linearext:morphform}. Choose splittings $s_W$ and $s_V$ of the extended linear spaces and represent $\tilde \phi$ by the pair $\phi : W \to V$ and $b_{\tilde \phi} \in \wedge^2 W^*$. Then we have the following diagram defining $\hat \phi : \wedge^* \hat V \to \wedge ^* \hat W$ and $\tilde \phi^\bullet : \mc S_V \to \mc S_W$:
$$\text{\xymatrix{
\wedge^* V^* \ar[d]_{\tilde \phi^* = e^{-b_{\tilde \phi}} \circ \phi^*} \ar[r]^{\hat{}} & \wedge^* \hat V \ar[d]^{\hat \phi} \ar@{=}[r] & \{s_V\} \times \wedge^* \hat V \ar@{=>}[r] & \mc S_V \ar[d]^{\tilde \phi^\bullet}\\
\wedge^* W^* \ar[r]_{\hat{}} & \wedge^* \hat W \ar@{=}[r] & \{s_W\} \times \wedge^* \hat W \ar@{=>}[r] & \mc S_W
}.}
$$

Now we show that the above definition does not depend on the choices made, i.e. the splittings $s_W$ and $s_V$. Choose another pair of splittings:
$$s_W'(X) = s_W(X) + \iota_X B_W \text{ and } s_V'(Y) = s_V(Y) + \iota_Y B_V,$$
then, for example, the induced identification of $\mc S_V \to \wedge^* V^*$ is given by post-composition with $e^{-B_V}$. Let $\rho_{\mc S} \in \mc S_V$ which corresponds to $\rho \in \wedge^* V^*$, and we follow the diagram:
$$\text{\xymatrix{
\rho \ar@{->}[r] & e^{-B_V} \rho \ar@{->}[r] & e^{-b_{\tilde \phi} - \phi^*B_V} \phi^*(\rho) \ar@{->}[r] & e^{-b_{\tilde \phi} - \phi^*B_V + B_W} \phi^* (\rho) \ar@{=}[r]
& e^{-b_{\tilde \phi}'} \phi^*(\rho)\\
}}$$
It follows that $\tilde \phi^\bullet(\rho_{\mc S})$ is well defined independent of the choice of splittings.
\qed

\section{Category of extended manifolds}\label{mfdext}

\subsection{Courant algebroids} \label{mfdext:courantalgebroids}
Following \cite{Kosmann5}, we have the following definition of a Courant algebroid:
\begin{defn}\label{mfdext:courant}
Let $E \to M$ be a vector bundle. A \emph{Loday bracket} $*$ on $\Gamma(E)$ is a $\RR$-bilinear map satisfying the Jacobi identity, i.e. for all $\mf X, \mf Y, \mf Z \in \Gamma(E)$,
\begin{equation}\label{mfdext:courantdefn1}
\mf X * (\mf Y * \mf Z) = (\mf X * \mf Y) * \mf Z + \mf Y * (\mf X * \mf Z).
\end{equation}
$E$ is a \emph{Courant algebroid} if it has a \emph{Loday bracket} $*$ and a non-degenerate symmetric pairing $\<,\>$ on the sections, with
an \emph{anchor map} $a : E \to TM$ which is a vector bundle homomorphism so that
\begin{eqnarray}
\label{mfdext:courantdefn2} a(\mf X) \<\mf Y, \mf Z\> & = & \<\mf X, \mf Y * \mf Z + \mf Z * \mf Y\> \\
\label{mfdext:courantdefn3} a(\mf X) \<\mf Y, \mf Z\> & = & \<\mf X * \mf Y, \mf Z\> + \<\mf Y, \mf X * \mf Z\>.
\end{eqnarray}
\end{defn}
The bracket in the definition is not skew-symmetric in general and the skew-symmetrization $[\mf X, \mf Y] = \mf X * \mf Y - \mf Y * \mf X$ is usually called the \emph{Courant bracket} of the Courant algebroid. The above definition is equivalent to the definition as given in, for example, \cite{Liu} or \cite{Gualtieri}. We rephrase the definition of generalized complex structure (\cite{Gualtieri, Hitchin}):
\begin{defn}\label{mfdext:extenddef}
An \emph{extended tangent bundle} $\TTC M$ is a Courant algebroid (cf. definition \ref{mfdext:courant}) which fits into the following extension:
$$0 \to T^*M \xto{a^*} \TTC M \xto{a} TM \to 0.$$
The extended tangent bundle $\TTC M$ is \emph{split} if the extension is split by some $s : TM \to \TTC M$, so that the image is isotropic.
An \emph{extended almost complex structure} $\JJC$ on $\TTC M$ is an almost complex structure on $\TTC M$ which is also orthogonal in the pairing $\<,\>$. Furthermore, $\JJC$ is \emph{integrable} and is called \emph{extended complex structure} if the $+i$-eigensubbundle $L$ of $\JJC$ is involutive with respect to either of the brackets $*$ or $[,]$.
\end{defn}

An example of extended tangent bundle is $\TT M = TM \dsum T^*M$ with the natural pairing and the bracket that Courant discovered \cite{Courant}, namely
$$[X+\xi, Y+\eta] = [X,Y] + \LLC_X \eta - \LLC_Y \xi -\frac{1}{2} d( \iota_X \eta - \iota_Y \xi ).$$
\begin{remark}
\rm{
An \emph{extended tangent bundle} $\TTC M$  is also known with the name of \emph{exact Courant algebroid}  \cite{Severa, Bursztyn}. The Loday bracket for an extended tangent bundle is also known as the \emph{Dorfman bracket}.
}\end{remark}

Now, the sequence
$$0 \to T^*M \xto{a^*} \TTC M \xto{a} TM \to 0$$
is always split in the sense of the definition \ref{mfdext:extenddef} and we have
$$\TTC M \simeq \TT M = TM \dsum T^*M : \mf X \mapsto X+\xi, \text{ with } X = a(\mf X), \xi = \mf X - s(X) \text{ for } \mf X \in \TTC M.$$
Such map $s$ is also called a \emph{connection} in \cite{Severa}. The choice of isotropic splitting determines a closed $3$-form
$$H(X, Y, Z) = 2\<s(X), [s(Y), s(Z)]\> = 2\<s(X), s(Y) * s(Z)\>.$$
Then the Courant algebroid $\TTC M$ is identified with the $H$-twisted Courant algebroid structure on $\TT M$. The space of isotropic splittings, which will be denoted $\mc I(M)$, is a torsor over $\Omega^2(M)$ and different choices give cohomologous $3$-forms. The action is given by $(B\circ s)(X) = s(X) + \iota_X B$. The class $[H] \in H^3(M, \RR)$ is the \emph{\v Severa class} of $\TTC M$ \cite{Severa}.
With this point of view, we regard the generalized tangent bundle $\TT M$ with $H$-twisted Courant bracket as the pair $(\TTC M, s)$ of extended tangent bundle with isotropic splitting (or equivalently, exact Courant algebroid with connection). In particular, we'll use $\TTC_0 M$ to denote $\TT M$ with standard Courant bracket (albeit splitting as well).
This point of view is justified in the sense that reduction by Hamiltonian action induces an extended manifold structure on the quotient which is not naturally split  \cite{Hu1}.

\subsection{Extended submanifolds}\label{mfdext:sub}
Let $i: F \subset M$ be a submanifold. Then the embedding defines an extended tangent bundle $\TTC F$ by $\TTC F = \Ann(\KKC)/\KKC$ where $\KKC = \Ann_{T^*M} TF \subset i^*T^*M$. We get the roof
\begin{equation}\label{mfdext:extendroof}\text{
\xymatrix{
&& \Ann(\KKC) \ar[lld]_\pi \ar@{^{(}->}[rrd]&&\\
\TTC F
&& && i^*\TTC M
}
}
\end{equation}
We note that \emph{a priori} the submanifold $F$ is not associated with an embedding $\TTC F \to i^*\TTC M$.
\begin{lemma}\label{mfdext:class}
Let $[H_M] \in H^3(M)$ be the \v Severa class of $\TTC M$, then $i^*[H_M] \in H^3(F)$ is the \v Severa class of $\TTC F$.
\end{lemma}
{\it Proof:}
Let $s_M : TM \to \TTC M$ be an isotropic splitting defining the twisting form $H_M$. The image of $TF$ under $s_M$ must lie in $\Ann(\KKC)$. Because $s_M(TF) \inter \KKC = \{0\}$, we see that $s_M$ descends to $s_F : TF \to \TTC F$, which is an isotropic splitting. The twisting form $H_F$ defined by $s_F$ is simply $i^*H_M$.
\qed

\begin{defn}\label{mfdext:extended}
An \emph{extended submanifold} $\mc F$ of $M$ with respect to the extended tangent bundle $\TTC M$
is a submanifold $i: F \subset M$ with an \emph{extended structure} $\mc E \subset \Ann(\KKC)$,
which is an involutive subbundle with respect to the bracket $[,]$, and which fits into the
following exact sequence
\begin{equation}\label{mfdext:extendedseq}
0 \to \mc C \to \mc E \xto{a} TF \to 0, \text{ for some subbundle } \mc C \subset i^*T^*M,
\end{equation}
where $a$ is the restriction of the anchor map for $M$. $\mc F$ is an \emph{isotropic} (resp.
\emph{non-degenerate}) extended submanifold if the restriction of $\<,\>$ to $\mc E$ vanishes
(resp. is non-degenerate).
\end{defn}
By definition, $\mc E$ descends to an involutive subbundle $\mc E_F = \mc E / (\KKC \inter \mc E)
\subset \TTC F$, which will be called the \emph{reduced structure}.
The set of extended structures on $F$ form a partially ordered set with ordering the inclusion of
$\mc E$.
We will consider in the following only \emph{maximal} extended structures in the partial order,
and will often drop the adverb \emph{maximally}.
We show in the following that maximally isotropic and non-degenerate structures correspond
respectively to the notion of \emph{generalized tangent bundle} in \cite{Gualtieri} and of
\emph{split submanifold} in \cite{Ben-Bassat}.


\begin{lemma}\label{mfdext:isotropic}
$F$ admits a (maximally) isotropic extended structure iff the twisting class of $\TTC F$ is $0$. For such $F$, the space $\mc I_F$ of (maximally) isotropic extended structures is a torsor over $\Omega^2_0(F)$. Upon a choice of splitting $s$ of $\TTC M$ defining the twisting form $H_s$, then there is $\tau_s \in \Omega^2(F)$ so that $d\tau_s + i^*H_s = 0$.
\end{lemma}
{\it Proof:} Suppose that $\mc F$ is a maximally isotropic
extended submanifold, then $\mc C = \KKC$. Choose a splitting $s$
of $\TTC M$ which defines a twisting form $H_s \in \Omega^3_0(M)$
and let $\mf X = s(X) + \xi \in \mc E$ for $X \in TF$ and $\xi \in
T^*M$, then we have
$$\<\mf X, \mf Y\> = \<s(X) + \xi, s(Y) + \eta\> = \iota_X\eta + \iota_Y \xi = 0 \text{ for } \mf X, \mf Y \in \mc E.$$
Define $\tau_s \in \Omega^2(F)$ by $\tau(X, Y) = \<\mf X, s(Y)\> = \iota_Y \xi$. We check that $\tau_s$ is well defined. For $\mf X' = s(X) + \xi' \in \mc E$, we have $\xi'-\xi \in \mc C = \KKC$ and $\iota_Y \xi' = \iota_Y \xi$. Let $s'$ be another splitting of $\TTC M$ so that $s'(X) = s(X)- \iota_XB$ for some $B \in \Omega^2(M)$, then the twisting form for $\TTC M$ becomes $H_{s'} = H_s - dB$ and we have $\tau_{s'} = \tau_s + i^* B$. It follows that $d\tau_s + i^*H_s$ is independent of the choice of splitting $s$.

We consider the reduced structure $\mc E_F = \mc E / \KKC$. The
sequence \eqref{mfdext:extendedseq} descends to $0 \to \mc E_F \to TF
\to 0$ in $\TTC F$. Since $\mc E$ is isotropic and involutive, we
see that $\mc E_F$ is too. It follows that the splitting $s_F : TF
\to \mc E_F \subset \TTC F$ gives a twisting form $H_F = \<s_F(X),
[s_F(Y), s_F(Z)]_F\>_F = 0$. In particular, the twisting class of
$\TTC F$ is $0$. Conversely, suppose that $\TTC F$ has twisting
class $0$ and let $s_F : TF \to \TTC F$ be a splitting defining
the twisting form $H_F = 0$. Consider the preimage $\mc E =
\pi^{-1}(s_F(TF))$, then the corresponding $\mc C$ in
\eqref{mfdext:extendedseq} is $\KKC$. For any $\mf X, \mf Y \in \mc
E$, we have $\pi([\mf X, \mf Y]) = [\pi(\mf X), \pi(\mf Y)]$ and
$\<\mf X, \mf Y\> = \<\pi(\mf X), \pi(\mf Y)\>$. Since $\mc E_F$
is isotropic and involutive, we see that $\mc E$ is isotropic and
involutive as well. We note that $\mc E$ as constructed is maximal
among the isotropic extended structures.

Let $\mc E'$ be another maximally isotropic extended structure on $F$, which induces the same splitting $\mc E_F$ of $\TTC F$. Let $\mf X \in \mc E$ and $\mf X' \in \mc E'$ so that $a(\mf X) = a(\mf X') = X \in TF$, then $\eta = \mf X' - \mf X \in T^*M$ and $\<\eta, TF\> = 0$. It follows that $\eta \in \KKC$ and $\mf X' \in \mc E$. Thus $\mc E$ is uniquely determined by the reduced structure $\mc E_F$. Since the splittings of $\TTC F$ defining twisting form $H_F = 0$ is a torsor over $\Omega^2_0(F)$, the same is true for $\mc I_F$.

We note that for $\mf X \in \mc E$, $s_F(X) = \pi(\mf X) = \pi(s(X) + \xi) = \pi(s(X)) + \iota_X \tau_s$ and $\LLC_X \eta - \iota_Y d\xi \in \KKC$ for $X, Y \in TF$ and $\eta, \xi \in \KKC$. We then compute
\begin{equation*}
\begin{split}
& \pi([\mf X, \mf Y]) = \pi(s([X, Y]) + \LLC_X \eta - \iota_Y d\xi + \iota_Y \iota_X H_s) \\
= & \pi(s([X, Y])) + \LLC_X \iota_Y \tau_s - \iota_Y d\iota_X \tau_s + \iota_Y \iota_X i^*H_s\\
= & s_F([X, Y]) + \iota_Y\iota_X (d\tau_s + i^*H_s).
\end{split}
\end{equation*}
It follows that $d\tau_s + i^*H_s = 0$ by the involutiveness of $s_F : TF \to \TTC F$.
\qed

\begin{lemma}\label{mfdext:nondeg}
A (maximally) non-degenerate extended structure on the submanifold $F$ is equivalent to an embedding of the induced extended tangent bundles (cf. \eqref{mfdext:extendroof}) $i_* : \TTC F \to \TTC M$ covering the embedding $i_*: TF \to TM$.
\end{lemma}
{\it Proof:}
Suppose that $\mc E$ is a maximally non-degenerate extended structure on $F$, then $\mc C \simeq T^*F$ via the restricted pairing $2\<,\>$. Since $\mc E$ is involutive, it is an extended tangent bundle over $F$ with the induced Courant algebroid structure. The projection $\pi|_{\mc E} : \mc E \to \TTC F$ is an isomorphism of extended tangent bundles and we obtain an embedding of extended tangent bundles $i_* = \pi|_{\mc E}^{-1}: \TTC F \to \TTC M$. The other direction is obvious.
%
\qed


\subsection{Product}
Let $(M, \TTC M)$ and $(N, \TTC N)$ be two smooth manifolds with extended
tangent bundles. The product $M \times N$ admits then natural extended
 tangent bundle $\TTC (M \times N)$ defined as follows. Let $\pi_i$ be the projection of $M \times N$ onto the $i$-th factor for $i = 1, 2$. Then as bundles, we have natural identifications
$$\TTC(M \times N) = \pi_1^* \TTC M \dsum \pi_2^* \TTC N \text{ and } T(M\times N) = \pi_1^* T M \dsum \pi_2^* T N.$$
The structure of Courant algebroid on $\TTC(M \times N)$ is then given by
declaring that the bracket and pairing all vanish between $\pi_1^* \TTC M$
 and $\pi_2^* \TTC N$. The axioms are easy to check and the extended tangent bundle $\TTC (M \times N)$ is defined.

\subsection{Reversion} \label{mfdext:rev} Similar to the linear case in \S \ref{linearext:reversion}, we have
\begin{defn}\label{mfdext:reversion}
Let $M$ be a smooth manifold with an extended tangent bundle $\TTC M$, on which the structure of a Courant algebroid is given by $(*, \<,\>, a)$, then the \emph{reversed} extended tangent bundle $-\TTC M$ is the same bundle with the structure of Courant algebroid given by $(*, -\<,\>, a)$. A bundle isomorphism $\tau :\TTC M \to \TTC M$ is called a \emph{standard reversion} if it gives an isomorphism of Courant algebroids $\tau : \TTC M \to - \TTC M$ and covers identity map on $TM$.
\end{defn}
Let $s$ be an isotropic splitting of $\TTC M$, which defines a
twisting form $H \in \Omega^3_0(M)$ as in \S
\ref{mfdext:courantalgebroids} and gives an identification $I: \TTC
M \to \TT M : \mf X = s(X) + \xi \mapsto X + \xi$, where $X =
a(\mf X)$ and $\xi = \mf X - s(X)$. Under this identification, we
have
$$I[\mf X, \mf Y] = [I(\mf X), I(\mf Y)]_H \text{ and } \<\mf X, \mf Y\> = \<X+ \xi, Y + \eta\>.$$
Since a reversion of $\TTC M$ only changes the sign of the pairing, under the same identification $I$, $-\TTC M$ has a $-H$-twisted Courant algebroid structure.
Let $\tau_s: \TT M \to \TT M : X + \xi \mapsto X - \xi$, then it is a reversion:
\begin{equation*}
\begin{split}
&[\tau_s(X + \xi), \tau_s(Y + \eta)]_H \\
= & [X, Y] - (\LLC_X \eta - \LLC_Y \xi - \frac{1}{2}d(\iota_X \eta - \iota_Y \xi) + \iota_Y \iota_X (-H)) \\
=& \tau_s([X+ \xi, Y+ \eta]_{-H}).
\end{split}
\end{equation*}
and $\<\tau_s(X+\xi), \tau_s(Y + \eta)\> = -\<X+\xi, Y + \eta\>$.
As in  \S \ref{linearext:reversion}, we may see that there is one-to-one correspondence between the set of reversions and the set of isotropic splittings of $\TTC M$. It follows that the set of reversions is a torsor over $\Omega^2(M)$.

Suppose that $\JJC$ is an extended complex structure on $\TTC M$, then it is also an extended complex structure on $- \TTC M$. When we want to be clear as to which extended tangent bundle we are looking at, we use $-\JJC$ to denote the one on $-\TTC M$.

\subsection{Morphisms}
Let $(M, \TTC M)$ and $(N, \TTC N)$ be two smooth manifolds with extended tangent bundles. Let $f : M \to N$ be a smooth map and $id \times f : M \to M \times N$ its graph.
On the product $M \times N$, we may endow another natural extended tangent bundle $\TTC(M \tilde \times N)$, which is given by the product of $(M, -\TTC M)$ and $(N, \TTC N)$.
An \emph{extended structure} on the map $f$ is defined to be an extended structure $\tilde {\mc E}$ on its graph $(id\times f)(M)$. More specifically, we have
$$0 \to \tilde {\mc C} \to \tilde {\mc E} \to TM \to 0,$$
where $\tilde {\mc E} \subset (id \times f)^*(\TTC (M \tilde \times N)) = -\TTC M \dsum f^*\TTC N$ is an involutive subbundle.
\begin{defn}\label{mfdext:morphdef}
The pair $(f, \tilde {\mc E})$ of map with extended structure is a \emph{morphism} if $\tilde {\mc E}$ is maximally isotropic, and we denote it as $\tilde f = (f, \tilde {\mc E}) : (M, \TTC M) \to (N, \TTC N)$.
\end{defn}
We note that the \v Severa class of $\TTC (M \tilde \times N)$ is given by $\pi_2^*[H_N] - \pi_1^*[H_M]$, then by lemma \ref{mfdext:isotropic}, we see that $f$ admits an extension into a morphism $\tilde f$ iff $f^*[H_N] = [H_M]$.
Similar to proposition \ref{linearext:morphform}, we have the following
\begin{prop}\label{mfdext:morphform}
Let $f : M \to N$ be a smooth map between extended manifolds, so
that $f^*[H_N] = [H_M]$. Choose splittings $s_M$ and $s_N$ of
$\TTC M$ and $\TTC N$ respectively, then any morphism $\tilde f :
(M, \TTC M) \to (N, \TTC N)$ extending $f$ is determined by an
element $b_{\tilde f} \in \Omega^2(M)$ so that $H_M = f^* H_N - d
b_{\tilde f}$. Let $B_M \in \Omega^2(M)$ and $B_N \in \Omega^2(N)$
and consider splittings $s'_M = s_M + \iota_\bullet B_M$ and $s'_N
= s_N + \iota_\bullet B_N$. Then the element $b'_{\tilde f} =
b_{\tilde f} - B_M + f^* B_N \in \Omega^2(M)$ determines the same
morphism $\tilde f$. \qed
\end{prop}
\begin{remark}\label{mfdext:nontrivial}
\rm{ From  proposition \ref{mfdext:morphform}, it is clear that with
a fixed splitting of $\TTC N$, a splitting of $\TTC M$ may be
chosen so that the $2$-form $b_{\tilde f}$ vanishes. For example,
let $F \subset M$ and $\TTC F$ defined by the diagram
\ref{mfdext:extendroof}, then the inclusion $i : F \to M$ can be
extended to a morphism $\tilde i$. In particular, choosing the
splittings as given in the lemma \ref{mfdext:class}, the two form
part $b_{\tilde i}$ vanishes. The same is not true vice versa,
because the map $f^* : \Omega^2(N) \to \Omega^2(M)$ is not
surjective in general. }\end{remark}
\begin{defn}\label{mfdext:composit}
Let $\tilde f : (M, \TTC M) \to (N, \TTC N)$ and $\tilde g : (N, \TTC N) \to (P, \TTC P)$ be two morphisms. Choose and fix splittings of the extended tangent bundles, then the \emph{composition} $\tilde h = \tilde g \circ \tilde f$ is given by $h = g\circ f$ and $b_{\tilde h} = b_{\tilde f} + f^*(b_{\tilde g}) \in \Omega^2(M)$.
\end{defn}

\subsection{The category} \label{mfdext:cat} We are now able to describe the category $\mc ESmth$ of extended manifolds. An object in $\mc ESmth$ is the pair $(M, \TTC M)$, a smooth manifold with extended tangent bundle. The \v Severa class of $\TTC M$ will be denoted $[H_M]$. A morphism $\tilde f : (M, \TTC M) \to (N, \TTC N)$ is given by an isotropic extended structure on the graph of a smooth map $f : M \to N$. We say that $\tilde f$ is an extension of $f$. The composition of morphisms is given by definition \ref{mfdext:composit}. We have the following fibration sequence:
$$\Omega^2_0(M) \to \Hom_{\mc ESmth}(M, N) \to C_\TTC^\infty(M, N)$$
where $C_\TTC^\infty(M, N)$ denote the space of smooth maps $f:M
\to N$ s.t $f^*[H_N] = [H_M]$. In particular, when $M = N$ as
extended manifolds, we may consider the group $\ms G_\TTC$ of
symmetries of $\TTC M$, consisting of invertible morphisms, which
fits in the following sequence:
$$0 \to \Omega^2_0(M) \to \ms G_\TTC \to \Diff_{[H_M]}(M) \to 1.$$
When the \v Severa class vanishes, we use the subscript $_0$ in the various notions.
 Then it is shown in \cite{Gualtieri} that $\ms G_0 = \Diff(M) \ltimes \Omega^2_0(M)$.
  In \S\ref{symmetry}, we will have more detailed discussion about the group
   $\ms G_\TTC$ and its Lie algebra. Note also that this is a category with involution, where the involution is given by the reversion $\TTC M \mapsto -\TTC M$.



\subsection{Spinors and cohomology}
Analogous to the linear algebra case, we have
\begin{defn}\label{spinor:abstract}The \emph{abstract spinor space}  for
$\TTC M$ is $$\mc S^\bullet(M) = \mc I(M) \times_{\Omega^2(M)}
\Omega^\bullet(M).$$
\end{defn}
Again, a choice of $s \in \mc I(M)$ gives an identification of
$\mc S^\bullet(M)$ to $\Omega^\bullet(M)$. Note that the space
$\mc S^\bullet(M)$ is no longer graded by $\ZZ$, instead, it only
remembers the induced $\ZZ_2$-grading. The wedge product of
differential forms makes $\mc S^\bullet(M)$ into a natural module
over $\Omega^\bullet(M)$.
\begin{lemma}\label{mfdext:diffgraded}
Let $s \in \mc I(M)$, $H$ the twisting form defined by $s$ and $\mc S^\bullet(M)$ identified with $\Omega^\bullet(M)$. Let $d_H = d - H \wedge$ on $\Omega^\bullet(M)$, then it defines a differential $d_\TTC$ on $\mc S^\bullet(M)$. The abstract spinor space $\mc S^\bullet(M)$ becomes a differential graded module over $\Omega^\bullet(M)$, via the wedge product of forms.
\end{lemma}
{\it Proof:} It is  easy to check that $d_H^2 = 0$. We will show
that it defines an operator on $\mc S^\bullet(M)$. Suppose that
$s' \in \mc I(M)$ is another splitting so that $s(X)- s'(X) =
\iota_X B$. Then via $s'$, the identification $\mc S^\bullet(M)
\to \Omega^\bullet(M)$ undergoes a $B$-transform and for $\rho \in
\Omega^\bullet(M)$ we have
\begin{equation*}
\begin{split}
d_{H'} (e^B \circ \rho) & = (d - (H-dB)\wedge) (e^{-B} \wedge \rho) \\
& = e^{-B}\wedge(-dB \wedge \rho + d\rho) - (H - dB)\wedge e^{-B} \wedge\rho \\
& = e^{-B} (d_H \rho) = e^B\circ(d_H \rho).
\end{split}
\end{equation*}
Let $\alpha \in \Omega^\bullet(M)$. Since $\wedge \alpha$ and the
$B$-transform on $\mc S^\bullet(M)$ commute, we see that $\wedge
\alpha$ is an operator on $\mc S^\bullet(M)$, which makes $\mc
S^\bullet(M)$ into a graded $\Omega^\bullet(M)$-module. For the
differential, we compute
$$d_H(\rho \wedge \alpha) = d(\rho \wedge \alpha) - H\wedge \rho \wedge \alpha = d_H\rho \wedge \alpha + (-1)^{|\rho|} \rho \wedge d\alpha.$$
\qed
\begin{defn}\label{mfdext:spindef}
The cohomology $H^\bullet(\TTC M)$ of $(\mc S^\bullet(M), d_\TTC)$ is the \emph{de Rham cohomology} of $\TTC M$.
\end{defn}
From the above lemma, we see that $H^\bullet(\TTC M)$ is a graded module over $H^\bullet(M)$. When we use $s \in \mc I(M)$ to identify $\mc S^\bullet(M)$ with $\Omega^\bullet(M)$, we see that the de Rham cohomology of $\TTC M$ is simply the {\bf $H$-twisted cohomology $H^\bullet(M; H)$ of $M$}. In particular, we see that $H^\bullet(M; H)$ depends only on the cohomology class $[H] \in H^3(M)$.

%
Let $\tilde f : (M, \TTC M) \to (N, \TTC N)$ be a morphism. Choose splittings of $\TTC M$ and $\TTC N$, so that $\tilde f$ is represented by $f: M \to N$ and $b_{\tilde f} \in \Omega^2(M)$. Define the \emph{pull-back} $\tilde f^\bullet : \mc S^\bullet(N) \to \mc S^\bullet(M)$ by
$$\tilde f^\bullet(\rho) = e^{-b_{\tilde f}} \wedge f^*(\rho).$$
Then it is independent of the choice of splitting as in the linear
case (proposition \ref{linearext:pullback}).
\begin{prop}\label{mfdext:pullback}
$\tilde f^\bullet$ is a chain homomorphism, and therefore it
induces  a homomorphism $\tilde f^\bullet : H^\bullet(\TTC N) \to
H^\bullet(\TTC M)$.
\end{prop}
{\it Proof:}
Choose splittings of $\TTC M$ and $\TTC N$ and let the twisting forms be $H_M$ and $H_N$ respectively. Then $H_M = f^*H_N - db_{\tilde f}$ and we compute
$$d_{H_M}(e^{-{b_{\tilde f}}} f^*(\rho)) =-db_{\tilde f}e^{-{b_{\tilde f}}} f^*(\rho) + e^{-{b_{\tilde f}}} d_{H_M} f^*(\rho) = e^{-{b_{\tilde f}}} d_{f^*H_N} f^*(\rho) = e^{-{b_{\tilde f}}} f^*(d_{H_N} \rho).$$
\qed

\subsection{Thom isomorphism} \label{mfdext:Thomiso}
Let $\pi: V \to M$ be an oriented real vector bundle of rank $k$.
In the classical situation, the Thom isomorphism is as follows.
Let $H^*_{cv}(V)$ denote the de Rham cohomology of forms with
compact support along the fibers, then in $H^k_{cv}(V)$, there is
a \emph{Thom class} $[\Theta]$, which is the unique class which
restricts to the orientation class $[\Theta_x] \in H^k_c(V_x)$.
The Thom isomorphism is therefore defined by wedging with
$[\Theta]$:
$$Th : H^*(M) \xto{\wedge [\Theta]} H^*_{cv}(V).$$

In the category $\mc ESmth$, a vector bundle $\tilde \pi : (V, \TTC V) \to (M, \TTC M)$ is a morphism extending a classical vector bundle $\pi$. It follows that the \v Severa class of $\TTC V$ is given by $[H_V] = \pi^*[H_M]$. Then splittings may be chosen so that the morphism $\tilde \pi$ is given by $\pi$ together with $b_{\tilde \pi} = 0$. It follows that $H_V = \pi^*H_M$. By the local to global principle we can show that the extended de Rham cohomology also has Thom isomorphism, via cupping with the Thom class. The ingredients in the argument are discussed in \S\ref{kirwan} for the equivariant case.

Let $U \subset M$ be a contractible open subset and $V_U \to U$ be the restriction
of $V$ on $U$. The Poincar\'e lemma then states that $H_U = H_M|_{U} = dB_U$
for some $B_U \in \Omega^2(U)$. Then it  follows that $\TTC h_U: H^*(U;H_U) \xto{\wedge [\Theta|_{U}]} H^*_{cv}(V_U;\pi^* H_U)$ is an isomorphism because we have
$$\TTC h_U = e^{-\pi^*B_U} \circ Th|_U \circ e^{B_U}.$$
 Applying the Mayer-Vietoris sequence to the opens sets $U, W$:
 $$\xymatrix{
   H^*(U \cap W;H_{U \cap W}) \ar[r] \ar[d]_{\cong}^{\wedge[\theta_{U \cap V}]} & H^*(U \cup W;H_{U \cup W})
  \ar[r] \ar[d]^{\wedge[\theta_{U \cup V}]} &
 H^*(U;H_{U })  \oplus H^*(W;H_{W })  \ar[d]_{\cong}^{\wedge[\theta_{U}] \oplus \wedge[\theta_{V}]}  \\
   H^*_{cv}(V_{U \cap W};\pi^*H_{U \cap W}) \ar[r]  & H^*_{cv}(V_{U \cup W};\pi^*H_{U \cup W}) \ar[r]  &
 H^*_{cv}(V_U;\pi^*H_{U })  \oplus H^*_{cv}(V_W;\pi^*H_{W })
   }$$
and with the use of  the five-lemma, we have that the Thom isomorphism on $U$ and $W$ implies the Thom isomorphism on $U \cup W$. Therefore
by induction on the open contractible sets that cover $M$, we have the Thom isomorphism in twisted cohomology, namely
$$H^\bullet(M;H)  \xto{\wedge [\Theta]} H^\bullet_{cv }(V;\pi^*H).$$ Then we we can conclude:
\begin{prop}\label{mfdext:thom}
Let $\tilde \pi$ be a vector bundle in the category $\mc ESmth$, then wedging with the Thom class induces an isomorphism
$$\TTC h : H^\bullet(\TTC M) \xto {\wedge [\Theta]} H^\bullet_{cv}(\TTC V).$$
\qed
\end{prop}

\section{Symmetries and actions}\label{symmetry}

\subsection{Generalized symmetries}\label{symmetry:generalized} For the triple $(\TT M, T^*M, \<,\>)$, where $\<,\>$ is the natural pairing, we have the group of symmetries $\mc G = \Diff(M) \ltimes \Omega^2(M)$. The action of $(\lambda, \alpha) \in \mc G$ on $\mf X = X + \xi \in \TT M$ is given by
$$(\lambda, \alpha) \circ \mf X = \lambda_*(X + \xi + \iota_X \alpha) = \lambda_*(\mf X + \iota_{a(\mf X)} \alpha),$$
where $\lambda_*$ on the forms is defined to be $(\lambda^{-1})^* = (\lambda^*)^{-1}$.
The composition law is
$$(\mu, \beta) \circ (\lambda, \alpha) = (\mu \circ \lambda, \lambda^*\beta + \alpha).$$
%
In the spirit of the present work, we think of $\mc G$ as the presentation of the \emph{group of generalized symmetries} of $\TTC M$ as a bundle with structures $(T^*M, \<,\>)$. As a corollary of the proposition \ref{mfdext:morphform}, the effect of choosing a different splitting gives an action of $B \in \Omega^2(M)$ on $\mc G$ by:
\begin{equation}\label{symmetry:Bgroup}
B \circ (\lambda, \alpha) = (\lambda, \alpha + \lambda^*B - B).
\end{equation}
It's easy to check that
$$B\circ \{(\mu, \beta) \circ (\lambda, \alpha)\} = \{B \circ (\mu, \beta)\} \circ \{B \circ (\lambda, \alpha)\}.$$
Thus, similar to the definition \ref{linearext:spinorabsp} for the abstract spinor space, we define the \emph{abstract group of symmetries} of $\TTC M$, without the Loday bracket $*$, as
$$\ms G = \mc I(M) \times_{\Omega^2(M)} \mc G.$$
A choice of a splitting $s \in \mc I(M)$ gives the identification of $\ms G$ with $\mc G$ by the section $\{s\} \times \mc G$ of the bundle $\mc I(M) \times \mc G$ over $\ms G$ with fiber $\Omega^2(M)$.

Let $\TTC M$ be an extended tangent bundle, whose \v Severa class is $[H_M]$ and $\ms G_\TTC$ the group of symmetries of $\TTC M$, as described in \S\ref{mfdext:cat}. Choose a splitting $s \in \mc I(M)$ which defines $H \in \Omega^3_0(M)$ with $[H] = [H_M]$.
The subgroup $\mc G_H \subset \mc G$ preserving the Courant bracket $[,]_H$ on $\TT M$ is given by
$$\mc G_H = \{(\lambda, \alpha) | \lambda^*H - H = d\alpha\}.$$
Then under the identification of $\{s\} \times \mc G$ with $\ms G$ we see that
$$\{s\} \times \mc G_H \text{ is identified with } \ms G_\TTC.$$
\begin{defn}\label{symmetry:genaction}
The action of a Lie group $G$ on $\TTC M$ is given by a Lie group homomorphism $\tilde \sigma : G \to \ms G_\TTC$. The action $\tilde \sigma$ is said to \emph{preserve the splitting} $s \in \mc I(M)$ if the image of $G$ lies in $\{s\} \times (\mc G_H \inter \Diff(M))$ under the identification given above.
\end{defn}

\begin{prop} \label{symmetry:preservesplit}
Let $G$ be a compact Lie group and $\tilde \sigma$ an action of it on $\TTC M$, then there exist $s \in \mc I(M)$ which is preserved by $\tilde \sigma$.
\end{prop}
{\it Proof:}
Let's start with a random splitting $s' \in \mc I(M)$, under which we obtain the homomorphism
$$\tilde \sigma': G \to \ms G \to \{s\} \times \mc G \simeq \Diff(M) \ltimes \Omega^2(M) : g \mapsto (\lambda_g, \alpha_g).$$
Let $H' \in \Omega^3(M)$ be the $3$-form determined by $s'$. From the equation \eqref{symmetry:Bgroup}, we are looking for a $B \in \Omega^2(M)$ so that
$$\alpha_g + \lambda_g^* B - B = 0 \text{ for all } g \in G.$$
Consider the action of $G$ on $\Omega^2(M)$ induced by $\tilde \sigma$:
$$g (B) := \alpha_g + \lambda_g^*B.$$
We check that
$$(gh) (B) = \alpha_{gh} + \lambda_{gh}^*B = \lambda_h^*\alpha_g + \alpha_h + (\lambda_g \lambda_h)^*B = g(h(B)),$$
i.e. it is indeed a well-defined action.
Because $G$ is compact, consider
$$B := \int_G h(0) d\mu(h) = \int_G \alpha_h d\mu(h),$$
where $\mu$ is the Harr measure normalized so that the volume of $G$ is $1$. Obviously we have for any $g \in G$
$$g(B) = \alpha_g + \lambda_g^* B = \int_G gh(0) d\mu(h) = B.$$
It follows that the splitting $s = s' + \iota_\bullet B$ is preserved by $\tilde \sigma$.
\qed

Now we go through the same process for the Lie algebras. The Lie algebra of $\mc G$ is $\mc X = \Gamma(TM) \dsum \Omega^2(M)$
with the Lie bracket
$$[(X, A), (Y, B)] = ([X, Y], \LLC_X B - \LLC_Y A).$$
The infinitesimal action of $(X, A)$ on $\mf Y \in \TT M$ is given by
$$(X, A) \circ \mf Y = - [X, Y] - \LLC_X \eta + \iota_Y A, \text{ where } \mf Y = Y+\eta \in \Gamma(\TT M).$$
We write down the \emph{$1$-parameter subgroup} generated by $(X, A)$:
$$e^{t(X, A)} := (\lambda_t, \alpha_t) \text{ where } \lambda_t = e^{t X} \text{ and } \alpha_t = \int_0^t \lambda_t^* A dt.$$

From \eqref{symmetry:Bgroup}, we see that the effect of choosing a different splitting gives an action of $B \in \Omega^2(M)$ on $\mc X$ by:
$$B \circ (X, A) = (X, A + \LLC_X B).$$
We also have
$$[B \circ (X, A), B\circ (Z, C)] = B\circ [(X, A), (Z, C)].$$
Thus, we define the Lie algebra $\ms X$ of the abstract group of symmetries of $\mc TTC M$, without the Loday bracket $*$, as:
$$\ms X = \mc I(M) \times_{\Omega^2(M)} \mc X.$$
Let $\mc X_H \subset \mc X$ be the Lie algebra of the group $\mc G_H$, then
\begin{defn}\label{symmetry:infinitesimal}
The \emph{abstract Lie algebra} $\ms X_\TTC$ of infinitesimal symmetries of $\TTC M$ is the Lie sub-algebra of $\ms X$ defined by $\{s\} \times \mc X_H$, where $H$ is the twisting form defined by $s \in \mc I(M)$.
\end{defn}
\noindent
From the definition, we see that $\ms X_\TTC$ is the Lie algebra of $\ms G_\TTC$.

In the following, we describe the relations among the Lie algebras $\ms X_\TTC$ where the structure of Courant algebroid varies. For this purpose, we fix a splitting $s \in \mc I(M)$ and work in terms of $\mc X$ and $\mc X_H$, etc.
When $H = 0$, we see that the Lie algebra of $\mc G_0$ is $\mc X_0 = \Gamma(TM) \dsum \Omega^2_0(M)$ with the standard bracket as given above.
%
Let the \emph{$H$-twisted Lie bracket} $[,]_H$ on $\mc X$ be given by
$$[(X, A), (Y, B)]_H = ([X, Y], \LLC_X B - \LLC_Y A + d\iota_Y\iota_X H).$$
Then the subspace $\mc X_0$ of $\mc X$ is  a Lie sub-algebra under this twisted bracket and the linear isomorphism
$$\psi_H : (\mc X, [,]) \to (\mc X, [,]_H) : (X, A) \mapsto (X, A + \iota_X H)$$
is in fact a Lie algebra isomorphism.
The following is straightforward (cf. \cite{Hu1}):
\begin{prop}\label{symmetry:twistbracket}
As Lie algebras, $(\mc X_H, [,]) = \psi_H^{-1} (\mc X_0, [,]_H)$.
\qed
\end{prop}
\noindent
\begin{defn}\label{symmetry:infsplitting}
Fix a splitting $s \in \mc I(M)$, which determines the $3$-form $H$. For $(X, A) \in \mc X_0$, the \emph{$H$-twisted infinitesimal action} is
$$(X, A) \circ ( Y + \eta)= - [X, Y] - \LLC_X \eta + \iota_Y (A - \iota_X H)$$
and the corresponding \emph{$H$-twisted $1$-parameter subgroup} is $e^{t(X, A-\iota_X H)}$. We say that the ($H$-twisted infinitesimal) action of $(X, A) \in \mc X_0$ \emph{preserves the splitting} $s$ if the corresponding ($H$-twisted) $1$-parameter subgroup is so, i.e. $A - \iota_X H = 0$.
\end{defn}

\subsection{Extended symmetries}  \label{extended_symmetries_kappa}
Let $\mf X \in \Gamma(\TTC M)$, then it defines an element in $\ms X_\TTC$. Choose a splitting $s \in \mc I(M)$, which defines the twisting form $H$, then we have:
$$\kappa : \Gamma(\TTC M) \to \mc X_H : \mf X = s(X) + \xi \mapsto \kappa(\mf X) = (X, d\xi - \iota_X H).$$
Following the identification of $\{s\}\times \mc X_H$ with $\ms X_\TTC$, we see that $\kappa$ is a well-defined map to $\ms X_\TTC$. The image of $\kappa$ in $\ms X_\TTC$ is the space of \emph{infinitesimal extended symmetries}, while the kernel is $\Omega^1_0(M)$, the closed $1$-forms. It's easy to see that the infinitesimal action of $\mf X$ via $\kappa$ on $\Gamma(\TTC M)$ is given by (the negative of) the (non-skew-symmetric) Loday bracket $- \mf X * \mf Y$. 
Let
$$(\lambda_t, \alpha_t) = e^{t(X, d\xi - \iota_X H)}$$
be the $1$-parameter subgroup of generalized symmetries generated by $\mf X$, then obviously the action of $(\lambda_t, \alpha_t)$ on $\TTC M$ does not depend on the choice of splitting either.
We write $e^{t\mf X} \subset \ms G_\TTC$ for the extended symmetries generated by $\mf X$ as above. We note that $e^{\Gamma(\TTC M)}$ form a subgroup $\ms E_\TTC$  of \emph{extended symmetries} of $\ms G_\TTC$. In fact, the Lie bracket on $\ms X_\TTC$ is compatible with either of the brackets on $\Gamma(\TTC M)$, i.e.
$$\kappa([\mf X, \mf Y]) = \kappa(\mf X * \mf Y) = [\kappa(\mf X), \kappa(\mf Y)].$$
\begin{defn} \label{mfdext:definitionextendedaction}
The action of the Lie algebra $\mf g$ on $\TTC M$ is given by a Lie algebra homomorphism $\tilde \sigma : \mf g \to \ms X_\TTC$. We say that $\mf g$ acts by \emph{extended symmetries} if $\tilde \sigma$ factors through $\kappa$, i.e. there is $\delta : \mf g \to \Gamma(\TTC M)$ so that $\tilde \sigma = \kappa\circ \delta$. Such an action is called \emph{isotropic} if the image of $\delta$ is isotropic with respect to $\<,\>$. The action $\tilde \sigma$ is \emph{integrable} to a $G$ action, if it is induced by a Lie group homomorphism $G \to \ms G_\TTC$. Then $G$ acts by \emph{extended symmetries}, or the $G$-action is \emph{isotropic}, if the corresponding $\mf g$-action is so.
\end{defn}

\subsection{Hamiltonian action} Let $(M, \JJC)$ be an \emph{extended complex manifold}, which is necessarily of even dimension $2n$. It is well known that the extended complex structure $\JJC$ induces a natural Poisson structure $\pi$ on $M$. Let $G$ be a connected Lie group acting on $M$ via a homomorphism $\tau: G \to \ms E_\TTC$, which is induced by a Lie algebra homomorphism $\tau_* : \mf g \to \ms X_\TTC$. The action $\tau$ is \emph{Hamiltonian} with \emph{moment map} $\mu: M \to \mf g^*$ if the geometrical action of $G$ is Hamiltonian with respect to the Poisson structure $\pi$, with equivariant moment map $\mu$, and the extended action is generated by $\JJC(d\mu)$.

%


\subsection{Action on spinors}\label{spinor}

The action of the Clifford bundle $Cl(\TTC M)$ on $\mc S^\bullet$ is defined
 via a choice of a splitting $s \in \mc I(M)$
 by $\iota_{\mf X} \rho = \iota_X \rho + \xi \wedge \rho$
 where $\xi = \mf X - s(X)$.
 By construction, $\iota_{\mf X} : \mc S^\bullet(M) \to \mc S^\bullet(M)$
 does not depend on the choice of splitting, nor on the Courant algebroid structure, and neither on $[H]$.

\begin{lemma}\label{equicoho:contract}
Choose a splitting $s \in \mc I(M)$ and define $\tilde\lambda_t : = e^{t\mf Y}$ and $$\LLC_{\mf X}^H \rho = -\left.\frac{d}{dt}\right|_{t=0} (\tilde \lambda_t \circ \rho) =  \LLC_X \rho + (d\xi - \iota_X H)\wedge \rho.$$ Then $\tilde\lambda_t$ and $\LLC_{\mf X}^H$ are well defined as operators on the space of spinors $\mc S^\bullet(M)$.
\end{lemma}
{\it Proof: }
We first recall from \cite{Gualtieri} that the effect of $B$-transform on $\rho \in \Omega^*(M)$ is given by $e^B \circ \rho = e^{-B}\wedge \rho$. Under the given splitting, the action of $e^{t\mf Y} = \tilde \lambda_t = (\lambda_t, \alpha_t)$ on $\Omega^*(M)$ are then given by $\tilde \lambda_t \circ \rho = \lambda_{t*} (e^{-\alpha_t} \rho)$. Suppose that $s' \in \mc I(M)$ is another splitting so that $s(X)- s'(X) = \iota_X B$, then under $s'$ we have
\begin{equation*}
\begin{split}
\tilde \lambda_t'  \circ (e^B\circ\rho) & = \lambda_{t*} (e^{-\alpha'_t} \wedge e^{-B} \wedge\rho) \\
& = \lambda_{t*} (e^{-(\alpha_t + \lambda_t^*B - B) - B} \wedge\rho)\\
& = e^{-B} \lambda_{t*}(e^{-\alpha_t} \wedge \rho) = e^B \circ (\tilde\lambda_t\circ \rho),
\end{split}
\end{equation*}
\begin{equation*}
\begin{split}
\left.\frac{d}{dt}\right|_{t=0} (\tilde \lambda_t \circ \rho) & = \left.\frac{d}{dt}\right|_{t=0} (e^{-\int_0^t\lambda_{r-t}^*(d\eta-\iota_Y H)dr} \wedge \lambda_{t*} \rho)\\
& = -\LLC_{Y}\rho + (-(d\eta-\iota_Y H) \wedge \rho) = -\LLC^H_{\mf Y}\rho.
\end{split}
\end{equation*}
The lemma follows.
\qed

\begin{defn}\label{equicoho:cartanop}
The \emph{system of Cartan operators} associated to $\TTC M$ on the space of spinors $\mc S^\bullet(M)$ is the collection $\LLC^\TTC_{\mf X}, \iota_{\mf X}$ and $d_\TTC$ defined as above.
\end{defn}

\begin{theorem}\label{equicoho:cartaneqs}
Let $\mf X * \mf Y$ be the Loday bracket for the extended tangent bundle $\TTC M$. Then we have:
\begin{eqnarray}
\label{equicoho:eqset1} \LLC^\TTC_{\mf X} \rho & = & (d_\TTC\iota_{\mf X} + \iota_{\mf X} d_\TTC)\rho \\
\label{equicoho:eqset2} [\LLC^\TTC_{\mf X}, \LLC^\TTC_{\mf Y}] \rho & = & \LLC^\TTC_{\mf X * \mf Y} \rho,\\
\label{equicoho:eqset3} (\iota_{\mf X} \iota_{\mf Y} + \iota_{\mf Y} \iota_{\mf X}) \rho & = & 2\<\mf X, \mf Y\> \rho, \\
\label{equicoho:eqset4} (\LLC^\TTC_{\mf X} \iota_{\mf Y} - \iota_{\mf Y} \LLC^\TTC_{\mf X})\rho & = & \iota_{\mf X * \mf Y} \rho, \\
\label{equicoho:eqset5} (d_\TTC \LLC^\TTC_{\mf X} - \LLC^\TTC_{\mf X}d_\TTC) \rho & = & 0, \\
\label{equicoho:eqset6} d_\TTC^2 \rho & = & 0.
\end{eqnarray}
\end{theorem}
{\it Proof: } \eqref{equicoho:eqset1} and \eqref{equicoho:eqset3}
follow directly from definition while \eqref{equicoho:eqset6} is
well-known fact. Let's first fix a splitting $s$ and consider in
$\TT M$ the action of generalized symmetries on the terms in the
above equations. Let $\tilde \lambda = (\lambda, \alpha) \in
\tilde\GGS$ be a generalized symmetry, then
\begin{equation}\label{equicoho:actions}
\tilde \lambda \circ H = \lambda_*(H-d\alpha), \text{ }
\tilde \lambda \circ \mf X  = \lambda_*(\mf X + \iota_X \alpha), \text{ and }
\tilde \lambda \circ \rho  = \lambda_*(e^{-\alpha}\rho).
\end{equation}
Thus we compute:
\begin{equation}\label{equicoho:btrans}
\begin{split}
d_{\tilde \lambda \circ H}(\tilde \lambda \circ \rho) & = \lambda_* (d - (H-d\alpha) \wedge)(e^{-\alpha} \wedge \rho) = \lambda_*(e^{-\alpha} \wedge (d - H\wedge) \rho) = \tilde \lambda \circ d_H \rho, \\
\iota_{\tilde \lambda \circ \mf X} (\tilde \lambda \circ \rho) & = \lambda_*(\iota_X(e^{-\alpha}\wedge \rho) + (\xi + \iota_X \alpha) \wedge e^{-\alpha} \wedge \rho) = \tilde \lambda \circ(\iota_{\mf X} \rho), \\
\LLC_{\tilde\lambda \circ \mf X}^{\tilde\lambda \circ H}(\tilde\lambda \circ \rho) & = \lambda_*(\LLC_X(e^{-\alpha} \wedge \rho) + (d\xi + d\iota_X \alpha - \iota_X(H-d\alpha)) \wedge e^{-\alpha} \wedge \rho) = \tilde\lambda\circ(\LLC_{\mf X}^H \rho). \\
\end{split}
\end{equation}

Let $\tilde\lambda_t = e^{t\mf Y}$ be the ($H$-twisted) $1$-parameter subgroup of generalized symmetries generated by $\mf Y \in \Gamma(\TTC M)$.
Let $H_t = \tilde\lambda_t \circ H$, then direct computation shows that $\frac{d}{dt} H_t = 0$ for all $t$ and it follows that $H_t = H$ for all $t$. Then \eqref{equicoho:btrans} provides
\begin{equation*}
d_\TTC(e^{t\mf Y} \circ \rho) = e^{t\mf Y} \circ (d_\TTC \rho), \text{ } \iota_{e^{t\mf Y} \circ \mf X}(e^{t\mf Y}\circ \rho) = e^{t\mf Y}\circ (\iota_{\mf X} \rho), \text{ and } \LLC^\TTC_{e^{t\mf Y} \circ \mf X}(e^{t\mf Y} \circ \rho) = e^{t\mf Y} \circ (\LLC^\TTC_{\mf X} \rho),
\end{equation*}
from which \eqref{equicoho:eqset2}, \eqref{equicoho:eqset4} and \eqref{equicoho:eqset5} follow.
\qed
\section{Extended equivariant cohomology}\label{equicoho}
\subsection{Cartan complex} Consider an
extended action $\tilde \sigma$ of a Lie algebra $\mf g$ on $\TTC
M$ as in definition \ref{mfdext:definitionextendedaction}. This
means that the map $\tilde \sigma$ factors through the map $\delta
:{\mf g} \to \Gamma(\TTC M)$ with $\tilde{\sigma} = \kappa \circ
\delta$. Then we can generalize the usual definition of the Cartan
complex to the extended case by considering the algebra $\mc S^\bullet (M)
\tensor \widehat{S}(\mf g^*)$ of formals series on $\mf g$ with values in $\mc S^\bullet(M)$.

The algebra $\widehat{S}(\mf g^*)$ is the $\mf a$-adic completion of the polynomial algebra $S(\mf g^*)$
where $\mf a$ is the ideal generated by all polynomials with zero constant term (we refer to chapter 10 of
\cite{AtiyahMacdonald} for the definition of the $\mf a$-adic completion). In practice what it means is that
if $\{u_j\}$ is a base for $\mf g^*$ and $S(\mf g^*)=\RR[u_1, u_2, \dots]$, then the $\mf a$-adic completion
is the algebra of formal series in the $u_j$'s:
$$\widehat{S}(\mf g^*) = \RR[[u_1, u_2, \dots]].$$
 The parity
of the elements are assigned according to the usual rule, i.e. the
forms are even or odd according to their degree while the formal series part
 $\widehat{S}(\mf g^*)$ is always even.

 There are various reasons why one needs to consider the algebra of formal series and not the polynomial algebra. One
 of them being the fact that the B-field transformations send forms to forms via multiplications with exponential maps. These
 exponential maps are in general not polynomial maps and therefore by considering only polynomials we would be reducing
 the change of coordinates transformations to a small group which is not the one we are interested in. Most of the algebraic properties
 of the Cartan model hold also for the completed model.

  We define
$$d_{\TTC, \delta} : \mc S^\bullet (M) \tensor \widehat{S}(\mf g^*) \to \mc S^\bullet
 (M) \tensor \widehat{S}(\mf g^*) : (d_{\TTC, \delta} \rho)(\tau)
 = d_\TTC \rho(\tau) - \iota_{\delta(\tau)}\rho(\tau) \text{ for } \tau \in \mf g,$$
then one may check that it's an odd operator and that
$$[d_{\TTC, \delta}^2(\rho)](\tau) = -\LLC^\TTC_{\delta(\tau)} \rho(\tau) + \<\delta(\tau), \delta(\tau)\>\rho(\tau).$$
Choosing dual basis $\{\tau_j\}$ and $\{u_j\}$ of the Lie algebra
$\mf g$ and its dual $\mf g^*$, we may rewrite the above equation
in coordinates:
$$d_{\TTC, \delta} \rho = d_\TTC \rho - \sum_j u_j \iota_{\mf X_j} \rho, \text{ and }$$
$$d_{\TTC, \delta}^2 \rho =
-\sum_{j} u_j \LLC_{\mf X_j}^\TTC \rho + \sum_{j,k} u_ju_k \<\mf
X_j, \mf X_k\> \rho,$$ where $\mf X_j := \delta(\tau_j)$. Notice
that if the extended action $\tilde \sigma$ is isotropic and if
$\LLC^\TTC_{\delta(\tau)} \rho =0$, then $d_{\TTC,
\delta}^2(\rho)=0$. Thus, following the definition of the Cartan
complex for equivariant cohomology, we propose
\begin{defn}\label{equicoho:eqdef} Let $\tilde \sigma$ be an
{\emph{isotropic}} extended action of a Lie algebra $\mf g$ on
$\TTC M$ that factors through $\delta : {\mf g} \to \Gamma(\TTC
M)$.
The {\bf{extended $\mf g$-equivariant Cartan complex}} is
\begin{equation}\label{equicoho:complex}
C^\bullet_{\mf g}(\TTC M; \delta) :=
\{\rho \in \mc S^\bullet(M) \tensor \widehat{S}(\mf g^*) | \LLC^\TTC_{\delta(\tau)} \rho(\tau) =
 0 \text{ for all } \tau \in \mf g\},
\end{equation}
with the odd differential $d_{\TTC, \delta}$. The cohomology
$H^\bullet_{G} (\TTC M; \delta)$ of the complex $C^\bullet_{\mf
g}(\TTC M; \delta)$ is the \emph{extended $\mf g$-equivariant de
Rham cohomology} of $\TTC M$ under the extended action $\tilde
\sigma$ defined by $\delta$. We will often drop the $\delta$ from
the notations.
\end{defn}

\noindent {\bf Notation:} The action $\tilde \sigma$  is always
integrable to a $G$-action for some Lie group $G$ with Lie algebra
$\mf g$. We will assume that such group $G$ is chosen so that no
confusion may arise. Therefore the extended $\mf g$-equivariant de
Rham cohomology will be denoted as $H^\bullet_G(\TTC M)$.

\begin{assumption}\label{equicoho:isotropyassum}
We assume that all the extended actions are isotropic.
\end{assumption}

\subsection{Module structure} Similar to lemma \ref{equicoho:contract},
we have the following for the equivariant case
\begin{lemma}\label{equicoho:module}
The complex $C^\bullet_{\mf g}(\TTC M)$ is a differential $\ZZ/2$-graded
module over the usual Cartan complex $C^\bullet_{\mf g}(M)$. It
follows that the cohomology $H^\bullet_{G}(\TTC M)$ is a $\ZZ/2$-graded
module over the usual equivariant cohomology $H^\bullet_{G}(M)$.
\end{lemma}
{\it Proof:} Choose a splitting of $\TTC M$ and let $H$ be the
3-form associated by the splitting. We only need to compute the
operators $\LLC^\TTC_{\delta(\tau)}$ and $d_{\TTC, \delta}$ for
$\rho \wedge \alpha$ where $\rho \in C^\bullet_{\mf g}(\TTC M)$
and $\alpha \in C^\bullet_{\mf g}(M)$:
\begin{equation*}
\begin{split}
\LLC^\TTC_{\delta(\tau)} (\rho \wedge \alpha) & = \LLC_{X_\tau}(\rho \wedge \alpha) + (d\xi_\tau - \iota_{X_\tau} H) \wedge(\rho \wedge \alpha)\\
& = \LLC^\TTC_{\delta(\tau)} \rho \wedge \alpha + \rho \wedge \LLC_{X_\tau} \alpha = 0,\\
d_{\TTC, \delta} (\rho \wedge \alpha) & = d(\rho \wedge \alpha) - H \wedge (\rho \wedge \alpha) - \sum_j u_j (\iota_{X_j} + \xi_j \wedge )(\rho \wedge \alpha) \\
& = d_{\TTC, \delta}\rho \wedge \alpha + (-1)^{|\rho|} \rho \wedge d_{\mf g} \alpha.
\end{split}
\end{equation*}
The last two equations imply the lemma. \qed

\subsection{Invariant function} Let $f: \mf g \to C^\infty(M)$ be a
$G$-equivariant linear map and denote $f_\tau \in C^\infty(M)$ the
image of $\tau$ under $f$. Then we may ``perturb'' the map $\delta
: \mf g \to \Gamma(\TTC M)$ by $df$:
$$\delta_f (\tau) = \delta(\tau) + df_\tau.$$
It's easy to check that $\tilde \sigma = \kappa \circ \delta_f$
and the extended action $\delta_f$ is again isotropic. We have
\begin{prop}\label{equicoho:equiB}
The two extended equivariant cohomologies $H^\bullet_G(M; \delta)$
 and \\ $H^\bullet_G(M; \delta_f)$, defined respectively from $\delta$
 and $\delta_f$, are isomorphic.
\end{prop}
{\it Proof:} Choose basis $\{\tau_j\}$ of $\mf g$ and dual basis $\{u_j\}$ of $\mf g^*$ and write $f_j = f_{\tau_j}$.
Let $b = \sum_j u_j f_j \in (C^\infty(M) \tensor \mf g^*)^G$ be the (ordinary) equivariant $2$-form representing $f$ and consider the \emph{equivariant $B$-transformation} $e^b \rho$:
$$d_{\TTC, \delta}(e^b \rho) = e^b(d_{\TTC, \delta} + d_G b \rho) = e^b(d_{\TTC, \delta} + \sum_j u_j df_j \rho) = e^b d_{\TTC, \delta_f} \rho.$$
By the module structure, we see that $e^b$ is a chain isomorphism between $(C^\bullet_G(\TTC M), d_{\TTC, \delta})$ and $(C^\bullet_G(\TTC M), d_{\TTC, \delta_f})$, and the proposition follows.
\qed

In fact, in the proof above, we may replace $b$ by any (ordinary)
equivariant $2$-form $$B_G = B +\sum_j u_j f_j \in \Omega^2(M)^G
\oplus (C^\infty(M) \otimes {\mf g}^*)^G  $$ and obtain:
\begin{lemma}
An equivariant $B_G$-transformation induces an isomorphism on the
extended $G$-equivariant cohomologies $H^\bullet_G(\TTC '
M,\delta_f) \cong H^\bullet_G(\TTC M, \delta)$ where the $3$-form defined by the same splitting $s$
for $\TTC M$ is $H$ and for $\TTC ' M $ is $H'= H -dB$. \qed
\end{lemma}

\subsection{Preserving a splitting} \label{equicoho:pure}
Let $G$ acts by extended action on $\TTC M$ and suppose that the action preserves a splitting $s \in \mc I(M)$, i.e.
$$d\xi_\tau - \iota_{X_\tau} H = 0 \text{ where } \delta(\tau) = s(X_\tau) + \xi_\tau \text{ for } \tau \in \mf g,$$
and $H$ is the $3$-form determined by $s$.
In such a case we have  that $\LLC_{\delta(\tau)} = \LLC_{X_\tau}$ and $d_{\TTC, G} = d_G - H_G \wedge$ where $d_G = d - \sum u_j \iota_{X_{\tau_j}}$ is the equivariant differential of the (ordinary) Cartan model, and $H_G := H + \sum_j u_j \xi_j$ is a closed equivariant 3-form also in the (ordinary) Cartan model. We then define
\begin{defn}
The cohomology of the equivariant Cartan complex
$$(\Omega^\bullet(M) \tensor \widehat{S}(\mf g^*))^{\mf g} = \{\rho \in \Omega^\bullet(M) \tensor \widehat{S}(\mf g^*) | \LLC_{X_\tau} \rho = 0 \text{ for all } \tau \in \mf g\} $$
with derivation $d_G - H_G \wedge$ will be called the {\bf
$H_G$-twisted equivariant cohomology of $M$}, and will be denoted
by $H^\bullet_{G}(M;H_G)$.
\end{defn}

 Thus we obtain
\begin{prop}\label{equicoho:twistident}
If the extended action $\tilde \sigma$ is isotropic and preserves a splitting, then $H^\bullet_G(\TTC M) = H^\bullet_{G}(M;H_G)$.
\end{prop}
{\it Proof.} As $$d_G H_G = dH +\sum_j u_j( d\xi_j - \iota_{X_{\tau_j}} H ) - \sum_{j,l} \<X_{\tau_j}+\xi_j, X_{\tau_l} + \xi_l \>u_ju_l$$
then $d_G H_G =0$ because $H$ is closed (i.e. $dH =0$), the action is pure, i.e. $d\xi_\tau - \iota_{X_\tau} H = 0$ and the action is isotropic, i.e.
$\<X_{\tau_j}+\xi_j, X_{\tau_l} + \xi_l \>=0$.
\qed

\begin{corollary}\label{equicoho:cpctlift}
If a compact Lie group $G$ acts by extended symmetries and is isotropic, then $H^\bullet_G(\TTC M)$ is isomorphic
to $H^\bullet_{G}(M;H_G)$ for some closed equivariant three form $H_G$.
\end{corollary}
{\it Proof:}
It follows from proposition \ref{symmetry:preservesplit}.
\qed

\subsection{Factor through a splitting}
In general, an ordinary closed $G$-invariant form $H$ (with
respect to the geometrical action) does not necessarily lift to an
equivariantly closed form $H_G$ with respect to $d_G$. From our
construction, when the action is, in certain sense, compatible
with the twisting, we can still obtain a cohomology encoding the
action that will be twisted by the (not necessarily invariant)
\emph{ordinary} closed $3$-form $H$. We note that when $H$ fails
to lift, according to corollary \ref{equicoho:cpctlift}, the $\mf g$-action
does not integrate to an action of a compact Lie group.

A $\mf g$-action
\emph{factors through a splitting $s$} if the map $\delta : \mf g
\to \Gamma(\TTC M)$ factors through $s$.
For such action, suppose further that $\iota_{X_\tau} \iota_{X_\omega} H = 0$ for all $\tau, \omega \in \mf g$, where $H$ is the twisting form defined by $s$. Then the extended $\mf g$-equivariant Cartan complex becomes
$$C^\bullet_{\mf g}(M; H) = \left\{\rho \in \Omega^\bullet(M) \tensor \widehat{S}(\mf g^*) | \LLC_{X_\tau} \rho - \iota_{X_\tau} H \wedge \rho = 0 \text{ for all } \tau \in \mf g\right\},$$
with differential $d_{G, H} = d_H - \sum_j u_j \iota_{X_j} = d_G - H\wedge$. The corresponding cohomology is denoted by $H_{H,G}^\bullet(M)$. We make the following definition:
\begin{defn}\label{equicoho:geotwist}
Let $M$ be an ordinary $G$-manifold and $H \in \Omega^3_0(M)$ so that $\iota_{X_\tau}\iota_{X_\omega} H = 0$ for all $\tau$, $\omega \in \mf g$, the \emph{$H$-twisted equivariant cohomology} of $M$ is defined to be $H_{H,G}^\bullet(M)$ as above.
\end{defn}

\subsection{Trivial action}\label{equicoho:trivial} We consider the
trivial action of $G$ on $\TTC M$, where $G$ can be taken as a compact Lie group.
This is an example of action preserving a splitting described in \S\ref{equicoho:pure}, while it is of independent
interest in the discussion of localization. In the ordinary case,
the trivial action of $G$ on $M$ gives the equivariant cohomology
$H_G^*(M;0) = H^*(M) \tensor \widehat{H}^*(BG)$, where $\widehat{H}^*(BG)=\widehat{S}(\mf g^*)^{\mf g}$ . Here, although the action of
$G$ on $\TTC M$ is trivial, the cohomology $H_{G}^\bullet(\TTC M)$
may be different from $H_G^*(M;0)$.

Consider a linear map $\delta : \mf g \to \Omega^1_0(M)$ where in
this case $\Omega^1_0(M) = \ker \kappa \subset \Gamma(\TTC M)$ as
defined in section \S \ref{extended_symmetries_kappa}. The induced
action of $G$ is then trivial on $\TTC M$ and the extended $\mf
g$-equivariant Cartan complex becomes:
$$C^\bullet_{\mf g}(\TTC M) = \Omega^\bullet(M) \tensor \widehat{S}(\mf g^*)
 \text{ and } d_{\TTC, \delta} \rho =
 d_\TTC \rho - \sum_j u_j \xi_j \wedge \rho,$$ where $\xi_j := \delta(\tau_j)$.
Choose a splitting $s$, which defines twisting form $H$, and we
rewrite the differential as:
$$d_{\TTC, \delta} \rho = d\rho - \left(H + \sum_j u_j \xi_j\right) \wedge \rho =
d\rho - H_G \wedge \rho,$$ where $H_G: = H + \sum_j u_j \xi_j$ is
seen as representing a cohomology class in $H_G^3(M)$ because
$d\xi_j =0$ and $X_{\tau_j}=0$. Thus, the extended equivariant
cohomology $H^\bullet_G(\TTC M)$ is the $H_G$-twisted equivariant
cohomology $H_{G}^\bullet(M;H_G)$.

\begin{corollary} \label{equicoho:cor_trivial}
If the 1-forms $\xi_j$ are exact (say when $H^1(M)=0$), then $H_G^\bullet(\TTC M)$ is isomorphic to the $H$-twisted cohomology of $M$ tensored with $\widehat{H}^\bullet(BG)$.
\end{corollary}
{\it Proof.}
There are functions $f_j$ over $M$ such that $df_j = \xi_j$. Then the B-field transform defined by $e^{-\sum_j u_j f_j}$ maps $H + \sum_j u_j \xi_j$ to $H$, then $H^\bullet_G(\TTC M) \cong H^\bullet (M; H) \tensor \widehat{H}^\bullet(BG)$.
\qed

\begin{example}  \label{example:trivial_action} {\rm
Take $M=S^1$ and $G=S^1$ with the trivial $G$ action in $M$, and
consider the extended action $\delta : \RR \to \Omega^1_0(M)$, $1
\mapsto d\theta$ with $\TTC M = TM \oplus T^*M$. The extended
equivariant Cartan complex becomes
$$C^\bullet_{\mf g}(\TTC M)=
 \Omega^\bullet(M) \otimes \RR [[u]] $$
with extended derivation $d_{\TTC, G } = d - u d\theta\wedge$.
%
An element $\alpha \in C^\bullet_{\mf g}(\TTC M)$ is of the
following form, where $f_i$ and $g_j$'s are functions over $S^1$:
$$\alpha = \sum_{i } f_i u^i + d\theta \sum_{j } g_j u^j \,\,\Rightarrow \,\,
d_{\TTC, G} \alpha = df_0 + \sum_{i > 0}(df_i - f_{i-1}d\theta) u^i .$$
Thus $\alpha$ is closed iff  $df_0 = 0$ and $df_i = f_{i-1}d\theta$ for all $i > 0$, which implies that $f_i=0$ for all $i$,
 i.e. the closed forms are $$\alpha =  d\theta\sum_{j } ^ng_j u^j \in C^{od}_{\mf g}(\TTC M).$$

 If we consider the form $\beta = \sum_{i } h_i u^i $ then $\alpha = d\beta$ is equivalent to the set of equations $dh_0=g_0 d\theta$,
 $dh_i = (g_i + h_{i-1})d\theta$, $i>0$ has a solution whenever $\int_{S^1}g_0 d \theta=0$. So we can conclude that

$$H^{ev}_G(\TTC S^1) = 0 \,\,\text{ and }\,\, H^{od}_G(\TTC S^1) =
 \RR.$$
This does not contradicts corollary \ref{equicoho:cor_trivial} as the $1$-form $d\theta$ is \emph{not} exact.
}\end{example}

\subsection{Circle bundle over surfaces}
Let $\pi: M \to \Sigma$ be an $S^1$-principle bundle over a closed
surface $\Sigma$ of genus $g$ and choose $H \in \Omega^3(M)$ be an
invariant volume form. We compute the cohomology group
$H^\bullet_G(M; kH)$, for $k \neq 0$, as defined in proposition
\ref{equicoho:geotwist}, which is a special case of the
equivariant cohomology $H_{G}^\bullet(\TTC M)$ for $\TTC M$ with
twisting class $[kH]$. We note that $H$ can not be lifted to an
equivariantly closed form in the usual Cartan model. Suppose that
$H_G = H + u \xi$ were such a lifting, then we compute $d_G H_G =
0 \iff d\xi = \iota_X H$ and $\iota_X \xi = 0$, where $X$ is the
infinitesimal action of $S^1$. It follows that $\xi =
\pi^*\lambda$ for some $\lambda \in \Omega^1(\Sigma)$. Then we
have $\int_M H = a\int_\Sigma d\lambda = 0$ for certain $a \neq
0$, which is contradiction as $H$ is a volume form.

The fact that $H$ fails to lift to an equivariantly closed form implies
that the geometrical action of $S^1$ cannot be lifted to an
isotropic extended action on $\TT M$.
It implies that the action under consideration is in fact an extended action of $\RR^1$.

We consider the complex:
$$C_G^\bullet(M; kH) = \{\rho \in \Omega^\bullet(M) \tensor
\widehat{S}(\mf g^*)| \LLC_X \rho - k\iota_X H \wedge \rho = 0\} \text{ and
} d_{G,H} = d -k H\wedge - u \iota_X.$$ Because $G = S^1$ is
abelian, we have
$$C_{\mf g}^\bullet(M; kH) = \{\rho \in \Omega^\bullet(M)| \LLC_X \rho - k\iota_X H \wedge \rho = 0\} \tensor S(\mf g^*).$$
Let $\theta$ be a connection form on $M$ and write $\rho \in \Omega^\bullet(M)$ as
$$\rho = \rho_0 + \theta \wedge \rho_1, \text{ where } \rho_i = a_i + \rho_i^1 + \rho_i^2 \text{ and } \iota_X \rho_i = 0,$$
with $a_i \in \Omega^0(M), \rho_i^j \in \Omega^j(M)$. Then $\rho
\in C_{\mf g}^\bullet(M; kH)$ iff
$$\LLC_X a_i = 0, \LLC_X \rho_i^1 = 0 \text{ and } \LLC_X\rho_i^2 = ka_i \iota_X H.$$
By $\iota_X \rho_i = 0$, we see that $\iota_X(d\rho_i^2 - a_i k H) = 0$, from which it follows that $d\rho_i^2 - a_i k H$ descends to $\Sigma$, i.e. $d\rho_i^2 = a_i k H$. Again, volume form argument implies that $a_i = 0$.
 In all, we have
$$\rho \in C_{\mf g}^\bullet(M; k H) \inter \Omega^\bullet(M) \iff \rho = \pi^*(\alpha_1+\alpha_2) + \theta\wedge \pi^*(\beta_1 + \beta_2).$$
where $\alpha_1, \beta_1 \in \Omega^1(\Sigma)$ and $\alpha_2, \beta_2 \in \Omega^2(\Sigma)$. Apply $d_{G, kH}$ we get
$$d_{G,kH} \rho =  \pi^* (d\alpha_1) - u \pi^*(\beta_1 + \beta_2) - d\theta_1 \wedge \pi^*(d\beta_1).$$
For a general form $\rho = \sum_j u^j \rho_j$ in the tensor product
we compute
$$d_{G,kH} \rho = - d\theta_1\wedge \pi^*\sum_{j=0} u^j d\beta_{1,j} - \pi^*\sum_{j=1} u^j(\beta_{1,j-1} + (\beta_{2, j-1} - d\alpha_{1,j})) + \pi^* d\alpha_{1,0}.$$
It follows that $d_{G,kH} \rho = 0$ is equivalent to
$$\beta_{1,j} = 0, \beta_{2, j} = d\alpha_{1,j+1} \text{ for all } j \vargeq 0 \text{ and } d\alpha_{1,0} = 0 \Rightarrow$$
$$\ker d_{G, kH} = \left\{\rho = d\theta_1 \wedge \pi^*\sum_{j = 0} u^j d\alpha_{1, j+1} + \pi^*\sum_{j = 0} u^j (\alpha_{1,j} + \alpha_{2,j}), \text{ with } d\alpha_{1,0} = 0 \right\}.$$
To conclude we find that the equivariant cohomology
$H^\bullet_G(M; kH)$ in the case of $k \neq 0$ is \emph{always}
the truncated de Rham cohomology of $\Sigma = M/G$:
$$H^\bullet_G(M; kH) \simeq H^\bullet \left\{0 \to \Omega^1(M/G) \xto{d} \Omega^2(M/G) \to 0\right\},$$
which, of course, maps to the usual de Rham cohomology of $\Sigma = M/G$.

For $k = 0$, the cohomology $H^\bullet_G(M; 0)$ is simply the usual equivariant cohomology, which is isomorphic to the de Rham cohomology of $\Sigma = M/G$.

\subsection{Non-free action on $S^3$} \label{section:non-free_action_example}
Let $S^3 \subset \CC^2$ be the unit sphere. We consider the standard coordinates
$z = (z_1, z_2) = (x_1 + i y_1, x_2 + iy_2) = (x_1, y_1, x_2, y_2)$ as well as
the polar coordinates
$(z_1, z_2) = r(e^{i\phi_1}\sin \lambda, e^{i\phi_2} \cos \lambda)$ on
$\CC^2$, where $r^2 = |z_1|^2 + |z_2|^2$, $\lambda \in [0, \frac{\pi}{2})$
and $\phi_j \in [0, 2\pi)$ for $j = 1, 2$.
Let $H = -\sin (2\lambda) d\lambda \wedge d\phi_1 \wedge d\phi_2$ and consider
the extended tangent bundle $\TTC M$ of $M = \CC^2\setminus \{(0,0)\}$ with
 \v Severa class $[H]$ (note that $[H] \neq 0$) with its corresponding splitting.
  Now, the embedding $i : S^3 \to M$ induces the extended structure $\TTC S^3$ with nontrivial \v Severa class and with the chosen splitting, we identify it as $\TT S^3$ with $i^*H$-twisted structure. We consider the action of $G = S^1$ on $S^3$ induced by rotating the first coordinate $z_1$:
$$\tilde \sigma : \RR^1 \to \Gamma(\TT S^3) : 1 \mapsto \mf X = \frac{\partial}{\partial \phi_1} - \cos^2\lambda d\phi_2,$$
which is pure with respect to the splitting.
Thus by proposition \ref{equicoho:pure}, the extended
$S^1$-equivariant cohomology $H_{S^1}^\bullet(\TTC S^3)$ is given
by the twisted equivariant cohomology $H^\bullet_{ S^1}(S^3; H_{S^1})$
with $H_{S^1} = H - u \cos^2 \lambda d\phi_2$.

\subsection{Calabi-Yau manifold}
We recall the notion of generalized Calabi-Yau manifold, which has a natural corresponding notion in the extended situation. For each maximally isotropic subbundle $L$ of $\TT_\CC M$, there is an associated spinor line bundle $U$, so that $L = \Ann U$ under Clifford multiplication. A generalized complex manifold $(M, \JJ)$ is called \emph{generalized Calabi-Yau} in the sense of \cite{Hitchin}, if there is $d$-closed non-vanishing section of $U$. Since Clifford multiplication is locally defined on the bundles, the definition of \emph{extended Calabi-Yau} manifold naturally extends. We show the following lemma, which is analogous to and a generalization of the corresponding one in symplectic geometry:
\begin{lemma}\label{equicoho:ham}
Let $(M, \JJC; \rho)$ be an extended Calabi-Yau manifold, i.e.
$\rho \in \Omega^\bullet(M)$ is non-vanishing with $\DDC\rho = 0$,
$\JJC$ corresponds to maximally isotropic subbundle $L \subset
\TTC_\CC M$, which anihilates $\rho$ under Clifford
multiplication. Suppose that there is Hamiltonian $G$-action on
$M$ with moment map $\mu$, then $\rho$ admits an equivariant
extension which is closed under $\DDC_G$.
\end{lemma}
{\it Proof: }
The main point is that $\<\JJC(d\mu), \JJC(d\mu)\> = \<d\mu, d\mu\> = 0$ which makes theorem \ref{equicoho:cartaneqs} suitable as an infinitesimal definition of a group action (in the sense of Cartan, see \S $2$ of Guillemin-Sternberg). Let $\mf X_j = \JJC(d\mu_j)$ where $\mu = \sum_j \mu_j u_j$, then by Hamiltonian-ness, we have $(\mf X_j + i d\mu_j) \cdot \rho = 0$, i.e. $\iota_{\mf X_j} \rho = -i d\mu_j \wedge \rho$. Let $\rho_G = e^{-i\mu} \rho$, then we compute that
\begin{equation*}
\DDC_G\rho_G =  \DDC(e^{-i\mu} \rho) - \sum_{j}u_j \iota_{\mf X_j} (e^{-i\mu} \rho) = 0.
\end{equation*}
\qed


\section{Properties of the extended equivariant cohomology}\label{kirwan}

In the first part we will show that the extended equivariant cohomology satisfies all the expected properties of a cohomology theory, and in the second we will show the Thom isomorphism when the action preserves a splitting, in particular, when the group is compact. We will not assume $G$ to be compact unless explicitly stated.

\subsection{Functoriality}
Let $(M, \TTC M)$ be an extended $G$-manifold, i.e., manifold $M$ with extended tangent bundle $\TTC M$ and an extended $G$-action. Let $(N, \TTC N)$ be another such manifold and consider the morphism $\tilde f = (f, \tilde {\mc E}) : (M, \TTC M) \to (N, \TTC N)$, where $\tilde {\mc E} \subset -\TTC M \dsum f^*\TTC N$ is the extended structure on $f$, then
\begin{defn}
$\tilde f$ is a \emph{extended $G$-equivariant morphism} if $f$ is $G$-equivariant in the usual sense and $\tilde {\mc E}$ is closed under the diagonal $G$-action on $-\TTC M \dsum f^*\TTC N$.
\end{defn}
\noindent In the alternative description of a pair $(f, b_{\tilde f})$ of map and $2$-form on $M$, we have
\begin{prop}\label{kirwan:equiv}
Choose and fix splittings of the extended tangent bundles $\TTC M$ and $\TTC N$, which determines a $b_{\tilde f} \in \Omega^2(M)$, then the morphism $\tilde f$ is equivariant iff
$$\xi_\tau = f^*\eta_\tau + \iota_{X_\tau} b_{\tilde f} \text{ for all } \tau \in \mf g,$$
where under the given splittings, $\mf X_\tau = \delta(\tau) = X_\tau + \xi_\tau$ for $\tau \in \mf g$, etc.
\end{prop}
{\it Proof:} For the action of $G$ to preserve $\tilde {\mc E}$,
we only need to show that the infinitesimal action of $\mf g$
preserves $\tilde {\mc E}$. The proposition then follows by
directly computing the elements of the form
$$\mf X_\tau *_M (X + f^*\eta + \iota_X b_{\tilde f}) + \mf Y_\tau *_N (f_*(X) + \eta), \text{ for } X \in TM \text{ and } \eta \in T^*N.$$
\qed

\begin{lemma} \label{kirwan:functorial}
Let $\tilde f$ be an equivariant morphism, then for $\rho \in \mc S^\bullet(N) \tensor \widehat{S}(\mf g^*)$,
\begin{enumerate}
\item $\LLC^\TTC_{\delta_M(\tau)} (\tilde f^\bullet(\rho)) = \tilde f^\bullet(\LLC^\TTC_{\delta_N(\tau)} \rho)$,
\item $d_{\TTC, \delta_M}(\tilde f^\bullet(\rho)) = \tilde f^\bullet(d_{\TTC, \delta_N} \rho)$.
\end{enumerate}
\end{lemma}
{\it Proof:} The map of linear spaces $\tilde f^\bullet: \mc
S^\bullet(N) \tensor \widehat{S}(\mf g^*) \to \mc S^\bullet(M) \tensor \widehat{S}(\mf
g^*)$ is induced by the map of  proposition
\ref{mfdext:pullback}. Let $\delta_M(\tau) = X_\tau + \xi_\tau$
and $\delta_N(\tau) = Y_\tau + \eta_\tau$, then both equations
follow from direct computation using the equations (under a fixed
choice of splittings of the extended tangent bundles):
$$Y_\tau = f_*(X_\tau), \xi_\tau = f^*\eta_\tau + \iota_{X_\tau} b_{\tilde f}, H_M = f^*H_N - db_{\tilde f} \text{ and } \tilde f^\bullet(\rho) = e^{-b_{\tilde f}} \wedge f^*\rho.$$
\qed

This lemma  \ref{kirwan:functorial} implies the functoriality of
the extended equivariant de Rham cohomology with respect to the
extended equivariant morphism $\tilde f$, so we have
\begin{corollary} \label{kirwan:pullback}
Let $\tilde f$ be an extended $G$ equivariant morphism $\tilde f :
(M,\TTC M) \to (N,\TTC N)$, then it induces a morphism
$$\tilde{f}^\bullet : H^\bullet_G(\TTC N, \delta_N) \to
H^\bullet_G(\TTC M, \delta_M).$$
\end{corollary}

\subsection{Mayer-Vietoris}
Let $U_1$ and $U_2$ be open subsets of $M$ with $U_1 \cup U_2 = M$ and let's denote by $U_{12}$ the intersection $U_1 \cap U_2$.  Then, it is clear that the sequence
$$0 \to (C^\bullet_{\mf g}(\TTC M), d_{\TTC,\delta_M}) \stackrel{j}{\to}
(C^\bullet_{\mf g}(\TTC U_1), d_{\TTC,\delta_{U_1}}) \oplus
(C^\bullet_{\mf g}(\TTC U_2), d_{\TTC, \delta_{U_2}})
\stackrel{k}{\to}  (C^\bullet_{\mf g}(\TTC U_{12}), d_{\TTC,
\delta_{U_{12}}}) \to 0$$ is exact, where $j(\omega) =
(\omega|_{U_1}, \omega|_{U_2})$ and $k(\alpha,\beta) = \alpha
|_{U_{12}} - \beta |_{U_{12}}$. Then it induces a long exact
sequence in cohomology

$$
\xymatrix{
H_G^\bullet(\TTC M) \ar[rr]^j & & H_{G}^\bullet(\TTC U_1) \oplus H_{G}^\bullet(\TTC U_2) \ar[rr]^k  & & H_{G}^\bullet(\TTC U_{12}) \ar[d] \\
H_{G}^{\bullet+1}(\TTC U_{12}) \ar[u] & &  H_{G}^{\bullet+1}(\TTC U_1) \oplus H_{G}^{\bullet+1}(\TTC U_2) \ar[ll]^k & &H_G^{\bullet+1}(\TTC M) \ar[ll]^j
}$$

\subsection{Exact sequence for a pair}

Suppose that there are extended $G$-actions $\tilde \sigma_F$ and
$\tilde \sigma_M$ on both $F$ and $M$. Let $\tilde i : (F, \TTC F)
\to (M, \TTC M)$ be a morphism extending the embedding $i: F
\subset M$, then $\tilde i$ is an \emph{equivariant embedding} if
the extended action $\tilde \sigma_F$ coincides with the action
induced from the roof \eqref{mfdext:extendroof}. Let $\mc K =
\Ann_{T^*M}TF$, then this is equivalent to the following:
$$\delta_M|_{F} : \mf g \to \Gamma(\Ann(\mc K)) \subset \Gamma(i^*\TTC M) \text{ and } \delta_F : \mf g \xto {\delta_M |_F} \Gamma(\Ann(\mc K)) \xto \pi \Gamma(\TTC F).$$
In particular, when we only consider the geometrical action, $F$
is an equivariant submanifold of $M$.
By corollary \ref{linearext:pullback}, we have the induced map of
equivariant cohomology $\tilde i^\bullet : H^\bullet_G(\TTC M) \to
H^\bullet_G(\TTC F)$, which exists at the chain level.

Now, performing the cone construction of $\tilde i^\bullet$, one
can define the relative complex as
$$C^\bullet_{\mf g}(\TTC M, \TTC F) : = C^\bullet_{\mf g}(\TTC M) \oplus C^{\bullet-1}_{\mf g}(\TTC F),
 \ \ \  d_{\TTC , {\delta}} (\omega, \theta) = (d_{\TTC M,\delta_M} \omega, \tilde i^\bullet \omega + d_{\TTC F,\delta_F} \theta)$$ that induces a short exact sequence of complexes
$$0 \to C^{\bullet-1}_{\mf g}(\TTC F) \to C^\bullet_{\mf g}(\TTC M, \TTC F) \to C^\bullet_{\mf g}(\TTC M) \to 0.$$

Then we get the long exact sequence in cohomology
$$
\xymatrix{
H_G^\bullet(\TTC M, \TTC F) \ar[rr] & & H_{G}^\bullet(\TTC M)  \ar[rr]^{\tilde i^\bullet}  & & H_{G}^\bullet(\TTC N) \ar[d] \\
H_{G}^{\bullet+1}(\TTC N) \ar[u] & &  H_{G}^{\bullet+1}(\TTC M)  \ar[ll]^{\tilde i^\bullet} & &H_G^{\bullet+1}(\TTC M, \TTC F) \ar[ll]
}$$
which is known as the exact sequence for a pair.

\subsection{Excision}
If we have the triple $A \subset Y \subset X$, the isomorphism
$H^\bullet_G(\TTC (X-A), \TTC (Y-A)) \cong H^\bullet_G(\TTC X,
\TTC Y)$ is obtained by using the  Mayer-Vietoris sequence for the
sets $X-A$ and $Y$, and the long exact sequences for the pairs
$(X,Y)$ and $(X-A, Y-A)$; this is an exercise in algebraic
topology.

\subsection{Thom isomorphism} \label{section:Thomiso}
This section is not completely satisfactory because we were not able to prove the Thom isomorphism
in the generality of the extended equivariant cohomology. Nevertheless we will show in what follows
the Thom isomorphisms for extended actions that preserves a splitting, for which the extended equivariant cohomology
is the twisted equivariant cohomology.

\subsubsection{Thom isomorphism for twisted equivariant cohomology}

Let $\pi : Z \to M$ be a $G$ equivariant real vector bundle of
rank $k$. Then integration along the fibers of the map $\pi$ on
the Cartan complexes
$$\pi_* :\left(\Omega^\bullet_{cv}(Z) \otimes S(\mf g^*) \right)^{\mf
g} \longrightarrow \left(\Omega^\bullet(M) \otimes S(\mf g^*)
\right)^{\mf g},$$ where $\Omega^\bullet_{cv}(Z)$ denotes
differential forms with vertical compact support, is a chain map
and induces an isomorphism in equivariant cohomology (see
\cite[Thm 10.6.1]{Guillemin})
$$\pi_*: H^{l}_{G,cv}(Z) \stackrel{\cong}{\longrightarrow} H^{l-k}_G(M).$$
This isomorphism is what is known as the Thom isomorphism. Its
inverse is obtained by wedging with the Thom class $\Theta \in
H^k_{G,cv}(Z) $
\begin{eqnarray*}
H^{l-k}_G(M) & \stackrel{\cong}{\to} & H^l_{G, cv} (Z) \\
\alpha & \mapsto & \Theta \wedge  \pi^* \alpha.
\end{eqnarray*}

Let $H_G$ be a closed equivariant three form on $M$ and let's
consider the $H_G$-twisted equivariant Cartan complex of $M$ and
the $\pi^*H_G$-twisted equivariant Cartan complex of $Z$ with
vertical compact support.

\begin{prop}
The map $\pi_*$ is chain map of twisted complexes
$$\pi_* : \left( \left(\Omega^\bullet_{cv}(Z) \otimes \widehat{S}(\mf g^*)
\right)^{\mf g}, d_G - \pi^*H_G \wedge \right) \longrightarrow
\left( \left(\Omega^\bullet(M) \otimes \widehat{S}(\mf
g^*)\right)^{\mf g} , d_G -H_G \wedge \right)$$ and it induces an
isomorphism of twisted equivariant cohomologies
$$ \pi_* : H_{G,cv}^*(Z, \pi^*H_G)
\stackrel{\cong}{\longrightarrow} H_G^*(M, H_G).$$
\end{prop}

{\it Proof.} Let's filter by degree the complexes
$$C_{\mf g, cv }(Z, \pi^*H_G) :=\left(\left(\Omega^\bullet_{cv}(Z) \otimes \widehat{S}(\mf g^*)
\right)^{\mf g}, d_G - \pi^*H_G \wedge \right)$$
$$C_{\mf g}(M, H_G):= \left(\left(\Omega^\bullet(M) \otimes \widehat{S}(\mf
g^*)\right)^{\mf g}, d_G -H_G\wedge \right)$$ with $\overline{F}^p
C_{\mf g}(M, H_G)$ the equivariant forms with degree $\geq p$ and
$F^p C_{\mf g, cv }(Z, \pi^*H_G)$ the equivariant forms of degree
$\geq p+k$.

The homomorphism $\pi_*$ is a chain map of twisted complexes
because
$$\pi_*(d_G \rho - \pi^*H_G \wedge \rho) = d_G \pi_* \rho -
\pi_*(\pi^*H_G \wedge \rho) = d_G \pi_* \rho - H_G \wedge
\pi_*\rho,$$ and moreover, it induces a morphism of filtered
differential graded modules
$$\pi_* : ( C_{\mf g, cv }(Z, \pi^*H_G), d_G - \pi^*H_G \wedge; F)
\to (C_{\mf g}(M, H_G), d_G -H_G\wedge; \overline{F}).$$

We have then that $\pi_*$ induces a homomorphism on the spectral
sequences associated to the filtrations $$\pi_* : E_k^{*,*} \to
\overline{E}_k^{*,*}$$ whose first level are the equivariant
differential forms
$$\pi_* : E_1^{*,*} =  \left(\Omega^\bullet_{cv}(Z) \otimes {S}(\mf
g^*)\right)^{\mf g} \to \overline{E}_1^{*,*} =
\left(\Omega^\bullet(M) \otimes {S}(\mf g^*)\right)^{\mf g}$$ and
and whose differential is equivariant derivative $\delta_1=d_G$.

Therefore the second level is the equivariant cohomology and
$\pi_*$ induces an isomorphism
$$\pi_* : E_2^{*,*}= H_{G,v}^\bullet (Z) \stackrel{\cong}{\to}
\overline{E}_2^{*,*} =H_G^\bullet(M).$$

Now, we also have that the twisted cohomology is complete
 with respect to the filtration, i.e. $$H^\bullet_{G, cv }(Z, \pi^*H_G)) =\lim_{\leftarrow p}
H^\bullet_{G, cv }(Z, \pi^*H_G)/F^p H^\bullet_{G, cv }(Z,
\pi^*H_G).$$ This last statement holds because of two facts: first
because the filtration by degree $F^pH^\bullet_{G, cv}(Z,
\pi^*H_G)$  is equivalent to the filtration $\mf a^pH^\bullet_{G,
cv}(Z, \pi^*H_G)$ given by the $\mf a$-dic topology,  where $\mf
a$ is the ideal of $S(\mf g^*)$ generated by polynomials with zero
constant term; and second because
 the twisted cohomology $H^\bullet_{G, cv}(Z,
\pi^*H_G)$ is complete with respect to the $\mf a$-dic completion,
as it is a finitely generated $\widehat{S}(\mf g^*)$-module.

 The facts that the twisted cohomologies are complete, that the filtrations are exhaustive
and weakly convergent (because the filtrations are by degree), and
that at the second level we have an isomorphism, imply by theorem
3.9 of \cite{McCleary} that $\pi_*$ induces an isomorphism of
twisted equivariant cohomologies
$$\pi_*: H^\bullet_{G, cv}(Z, \pi^*H_G) \stackrel{\cong}{\to}
H_G^\bullet(M, H_G).$$ \qed

By the same argument as in the untwisted case, the inverse map of
$\pi_*$ is given by wedging with the equivariant Thom form
$\Theta$. We can conclude that we have an isomorphism
$$Th: H^\bullet_{G}(M, H_G) \xto{\wedge [\Theta]} H^\bullet_{ G, cv}(Z, \pi^*H_G).$$

\subsubsection{Thom isomorphism for pure extended actions}

Let $\tilde \pi: (Z, \TTC Z) \to (M, \TTC M)$ be a vector bundle
in the category $\mc ESmth$, where a Lie group $G$ acts on $(M,
\TTC M)$ by generalized (resp. extended) symmetries. The bundle
$\tilde \pi$ is an equivariant bundle if there is a generalized
(resp. extended) $G$-action on $(Z, \TTC Z)$ lifting the one on
$M$. It means that the action on $Z$ is fiberwise linear and the
action on $M$ is induced from that of $Z$ by restricting to the
$0$-section.

When the $G$ action on $\TTC M$ preserves a splitting,
we see that $H^\bullet_G(\TTC M) \cong H^*_{G}(M;H_G)$ and $H^\bullet_{G,cv}(\TTC Z) \cong H^*_{G,cv}(Z;\pi^*H_G)$. As
$$Th : H^\bullet_{ G}(M;H_G) \stackrel{\cong}{\to} H^\bullet_{G,cv}(Z;\pi^*H_G)$$
we have that it induces an isomorphism of extended equivariant cohomologies
$$Th: H^\bullet_G(\TTC M) \stackrel{\cong}{\to} H^\bullet_{G, cv}(\TTC Z).$$

\section{Localization}\label{local}
\subsection{Fixed point set of generalized action}
Let's fix the splitting $s \in \mc I(M)$ that is preserved by the action of the (not necessarily compact) Lie group $G$.
Then by definition, the image of the homomorphism
$$\tilde \sigma : G \to \ms G_\TTC \simeq \{s\} \times \mc G_H$$
lies in $\mc G_H \inter \Diff(M)$.
Let $x \in F \subset M$ be a fixed point where $F$ is a component of fixed point set. Consider the induced representation of $G$ on $\TT _xM$. The action preserves the splitting of $\TT_x M = T_xM \dsum T_x^*M$. Furthermore, the representations $T_xM$ and $T_x^* M$ split as following:
$$T_xM = T_xF \dsum N_x \text{ and } T^*_xM = T^*_xF \dsum N^*_x,$$
where we have $T^*_xF = N_x^\perp$ and $N^*_x =T_xF ^\perp$ with respect to the natural
 pairing between $TM$ and $T^*M$. This can be seen in what follows. It is obvious that $T_x M$ splits as such, where the $T_xF$ component is simply the trivial sub-representation, while $N_x$ is the non-trivial part. For $T^*_xM$, we choose dual basis $\{v_i\}$ and $\{u_i\}$ of $T_xM$ and $T^*_xM$ respectively, so that $T_xF = \Span(v_{i = 1, \ldots, k})$ and $N_x = \Span(v_{i = k+1, \ldots, n})$. Let $g \in G$ and we compute
$$\left(g_{ij}\right) = \left(\<g\circ u_i, v_j\>\right) = \left(\<u_i, g\circ v_j\>\right) = \left( \begin{matrix}I_{k\times k} & 0 \\ 0 & * \end{matrix}\right)$$
In particular, it follows that $g\circ u_i = u_i$ for $i = 1, \ldots, k$. Let $T^*_xF = \Span(u_{i = 1, \ldots, k})$ and $N^*_x = \Span(u_{i = k+1, \ldots, n})$, then they are sub-representations of $T^*_xM$ and $N^*_x$ is the non-trivial part. Thus the representation $\TT_x M = \TT_x F \dsum \NN_x$ naturally splits into trivial and non-trivial components.

Let $\TT F = \union_{x \in F} \TT_x F$, then we show that
\begin{lemma}\label{equicoho:fixedclosed}
$\Gamma(\TT F)$ is closed under the Courant bracket.
\end{lemma}
{\it Proof:}
We consider the induced homomorphism of Lie group:
$$\tilde \sigma : G \to \ms G_\TTC \simeq \{s\} \times \mc G_H : g \mapsto (\lambda_g, \alpha_g).$$
Because $\tilde \sigma$ fixes the splitting $s$, we have $\alpha_g = 0$ and $\lambda_{g}^* H = H$ for all $g \in G$.
Then by definition,
$$\mf Y = Y + \eta \in \Gamma(\TT F)\iff \lambda_{g*} Y = Y \text{ and } \lambda_{g*} \eta = \eta \text{ for all } g \in G.$$
Let $\mf Z \in \Gamma(\TT F)$, straight forward computation gives $\lambda_{g*} [\mf Y, \mf Z]_H = [\mf Y, \mf Z]_H$.
\qed

By lemma \ref{mfdext:nondeg}, we see that $\TT F$ is isomorphic to
the induced extended tangent bundle $\TTC F$ in
\eqref{mfdext:extendroof} as a Courant algebroid. Thus the \v Severa
class of $\TTC F$ is given by $[H_F] = [i^* H_M] \in H^3(F)$.

\begin{corollary}
Let $\tilde \sigma : G \to \ms G$ be a proper extended action and $F$ a
fixed point set component of the geometrical action $\sigma$.
Then the induced action on $\TTC F$ is trivial, as described in \S\ref{equicoho:trivial}. \qed
\end{corollary}

\subsection{Localization in twisted equivariant cohomology}

In this section we show the localization theorem in twisted
equivariant cohomology following \cite{AtiyahBott}.
Here we restrict to the case where $G$ is compact. Then proposition
\ref{symmetry:preservesplit} implies that we may fix a splitting $s \in \mc I(M)$ and assume that the
$G$-action preserves $s$.

As we have seen in proposition \ref{equicoho:twistident} the
cohomology of $H^\bullet_G(\TTC M)$ could be calculated using the
twisted equivariant de Rham cohomology
$$(\Omega^\bullet(M) \tensor \widehat{S}(\mf g^*))^G = \{\rho \in \Omega^\bullet(M)
 \tensor \widehat{S}(\mf g^*) | \LLC_{\delta(\tau)} \rho = 0 \text{ for all } \tau \in \mf g\}
  \text{ and } \DDC_G = d_G - \alpha\wedge,$$
twisted by $\alpha: = H_G=H + \sum_j u_j \xi_{\tau_j}$ whenever the
action is pure, i.e. $d\xi_{\tau_j} - \iota_{X_{\tau_j}}H=0$
Recall that when the action is pure $\alpha$ defines a cohomology
class $[\alpha] \in H^3_G(M)$ and the twisted equivariant de Rham
cohomology is denoted by $H^\bullet_{ G}(M;\alpha)$.

As $H^\bullet_{G}(M;\alpha)$ is a module over $H^\bullet_G(M)$
(see lemma  \ref{equicoho:module}) then $H^\bullet_{G}(M;\alpha)$
becomes a module over $H^\bullet(BG) = S(\mf g^*)^{\mf g}$. If $i : F \to M$ is the
inclusion of the fixed point set of the geometrical action of $G$,
we will show that the pullback $i^*$ and the pushout $i_*$ in
extended equivariant cohomology are inverses of each other after
inverting the equivariant Euler class of the normal bundle of $F$.
For this we will mimic the proof of Atiyah and Bott of the
localization theorem in equivariant cohomology \cite{AtiyahBott,
Audin}.

For the sake of simplicity we will focus on the case that $G$ is a
torus and we will make use of complex coefficients. Having setup
the hypothesis we can start.

Recall that for $G$ a torus we have that $H^\bullet(BG) = S(\mf
g^*) = \mathbb{C}[u_1, \dots u_n]$ a polynomial ring in $n$
variables. The support of a $H^\bullet(BG)$-module $A$ is,
$supp(A) = \bigcap_{\{f| f\cdot A=0\}} V_f$ where $V_f=\{ X \in
\mf g | f(X)=0 \}$, then
\begin{lemma}
$supp( H^\bullet_{ G}(M;\alpha)) \subset supp (H^\bullet_{ G}(M))$
\end{lemma}
{\it Proof.} If $f\cdot H^\bullet_{ G}(M) = 0$, by the definition
of the $H^\bullet(BG)$-module structure in $H^\bullet_{
G}(M;\alpha)$, one has that $f \cdot H^\bullet_{G}(M;\alpha)=0$. The
inclusion follows. \qed

Let $F=M^G$ be the fixed point set of the $G$-action, then from
\cite[Prop 5.2.5]{Audin} we know that $supp ( H^\bullet_{G}(M-F))
\subset \bigcup_{H} \mf h $ where $H$ describes the finite set of
proper stabilizers of points in $M-F$, and $\mf h$ is the Lie
algebra of $H$. Then we have that $supp (
H^\bullet_{G}(M-F;\alpha)) \subset \bigcup_{H} \mf h $, and
therefore $H^\bullet_{G}(M-F;\alpha)$ is a torsion
$H^\bullet(BG)$-module.

\begin{lemma}
The kernel and the cokernel of the map $i^*:
H^\bullet_{G}(M;\alpha) \to H^\bullet_{G}(F;i^*\alpha)$ have
support contained in $\bigcup_{H} \mf h$ where the $H$ runs over
the stabilizers $\neq G$ of points in $M$.
\end{lemma}

{\bf Proof.} Let $U$ be an equivariant tubular neighborhood of
$F$. We know that $supp ( H^\bullet_{G}(M-U;\alpha)) \subset
\bigcup_{H} \mf h $ and also $supp (
H^\bullet_{G}(\partial(M-U);\alpha)) \subset \bigcup_{H} \mf h $.
Using the long exact sequence for the pair $(M-U, \partial(M-U))$
we conclude that
$$supp ( H^\bullet_{G}(M-U,\partial(M-U);\alpha))
\subset \bigcup_{H} \mf h .$$
 Let $V$ be another equivariant
tubular neighborhood containing $U$, such that $V-U \sim
\partial(M-U) =\partial U$. Then by excision
$H^\bullet_{G}(M,F;\alpha)\cong H^\bullet_{G}(M,V;\alpha) \cong
H^\bullet_{G}(M-U,\partial(M-U);\alpha)$, so in particular
$H^\bullet_{G}(M,F;\alpha,) \subset \bigcup_{H} \mf h $. From the
long exact sequence for the pair $(M,F)$,
$$H^\bullet_{G}(M,F;\alpha,) \to H^\bullet_{G}(M;\alpha)
\stackrel{i^*}{\to} H^\bullet_{G}(F;i^*\alpha,) \to
H^\bullet_{G}(M,F;\alpha)$$ the lemma follows. \qed

The Thom Isomorphism in twisted cohomology is obtained by wedging
with the Thom class (see \S \ref{section:Thomiso}). So for $\pi:
\nu \to F$ the equivariant normal bundle of $F$ in $M$ of rank
$d=dim(M) - dim(F)$ (seen as a tubular neighborhood $q : \nu \to
M$), and $\Theta \in H^d_{cv,G}(\nu)$ the Thom class, the Thom
isomorphism is
\begin{eqnarray*}
th : H^\bullet_{G}(Fi^*\alpha) &\to &
H^\bullet_{cv,G}(\nu;\pi^*i^*\alpha)\\
a &\mapsto & \pi^*(a) \wedge \Theta.
\end{eqnarray*}
But we need to land in $H^\bullet_{cv,G}(\nu;q^*\alpha)$. Then we
use the fact that $\pi^*$ is in isomorphism in equivariant
cohomology and therefore there is an equivariant form $\sigma$ on
$\nu$ such that $q^*\alpha = \pi^* i^* \alpha - d_G \sigma$. This
gives us the isomorphism $H^\bullet_{cv,G}(\nu;\pi^*i^*\alpha)
\stackrel{e^{-\sigma}}{\to} H^\bullet_{cv,G}(\nu;q^*\alpha)$, that
precomposed with the Thom map is what we are going to call (by
abuse of notation) the Thom isomorphism $Th : = e^{-\sigma} \circ
th$.

The pushforward map $i_* : H^\bullet_{G}(F;i^*\alpha) \to
H^\bullet_{G}(M;\alpha)$ is defined as the composition of the maps
$$H^\bullet_{G}(F;i^*\alpha) \stackrel{Th}{\to}
H^\bullet_{cv,G}(\nu;q^*\alpha)\cong
H^\bullet_{G}(M,M-F;\alpha) \to  H^\bullet_{G}(M;\alpha),$$ and
recall that the equivariant Euler class $e_G(\nu)$ is defined as
the element in $H^d_{G}(F)$ such that $\pi^*
e_G(\nu) = j^* \Theta$, where $j^*: H^\bullet_{cv,G}(\nu) \to H^\bullet_G(\nu)$ is the natural homomorphism

\begin{lemma} \label{pullbackpushforward}
The composition $i^*i_*$ is equivalent to multiplying by
$e_G(\nu)$ (using the $H^\bullet_G(F)$-module structure), i.e.
 for $a \in H^\bullet_{G}(F;i^*\alpha)$, we have $i^*i_*(a) = a
 \wedge e_G(\nu)$.
\end{lemma}

{\it Proof.} Consider the commutative diagram
$$
\xymatrix{H^\bullet_{G}(F;i^*\alpha) \ar[r]_{Th} \ar
@/^1pc/[rr]^{i_*} & H^\bullet_{cv,G}(\nu;q^*\alpha) \ar[r]
\ar[d]^{j^*} &
H^\bullet_{G}(M;\alpha) \ar[d]^{i^*}\\
 & H^\bullet_{G}(\nu;q^*\alpha) &
H^\bullet_{G}(F;i^*\alpha). \ar[l]^{\pi^*} }$$ Then from the left
side  $j^* Th(a) = \pi^*a \wedge j^*\Theta= \pi^*(a \wedge e_G(\nu))$,
and from the right hand side $j^* Th(a) = \pi^*(i^*i_*(a))$, and
as $\pi^*$ is an isomorphism, one obtains that $i^*i_*(a) = a
 \wedge e_G(\nu)$. \qed

\begin{lemma}
The kernel and the cokernel of the map
$i_*:H^\bullet_{G}(F;i^*\alpha) \to H^\bullet_{G}(M;\alpha)  $
have support contained in $\bigcup_{H} \mf h$ where the $H$ runs
over the stabilizers $\neq G$ of points in $M$.
\end{lemma}

{\it Proof.} Consider the cohomology exact sequence of the pair
$(M,M-F)$
$$\xymatrix{
H^\bullet_{G}(M-F;\alpha) \ar[r] & H^\bullet_{G}(M,M-F;\alpha)
\ar[d]^{\cong} \ar[r] & H^\bullet_{G}(M;\alpha) \ar[r]  &
H^\bullet_{G}(M-F;\alpha)\\
& H^\bullet_{G}(F;i^*\alpha) \ar[ur]_{i_*} & & }$$
and as $supp(H^\bullet_{G}(M-F;\alpha)) \subset \bigcup_H \mf h$
the lemma follows.
\qed

For $Z \subset F$ a connected component of the fixed point set
there exist a polynomial $f_Z \in H^\bullet(BG) $ such that the
Euler class $e_G(\nu_Z)$ of $\nu|_Z$ is invertible in
$(H^\bullet_{G}(F))_{(f_Z)}$, the localization of
$H^\bullet_{G}(F)$ in the ideal generated by $f_Z$ as a
$H^\bullet(BG)$-module. Therefore for $f=\prod_{Z\subset F} f_Z$
then $e_G(\nu)$ is invertible in $(H^\bullet_{G}(F))_{(f)}$ and
moreover the kernel of $f$ is contained in $\bigcup_{H \neq G} \mf
h$.

Now, define the homomorphism $Q:(H^\bullet_{G}(M;\alpha))_{(f)} \to
(H^\bullet_{G}(F;i^*\alpha))_{(f)}$ by $Q(a): = \sum_{Z \subset F}
i^*_Z(a) \wedge e_G(\nu_Z)^{-1}$ where $e_G(\nu_Z)^{-1} \in
H^\bullet_{G}(F)_{(f)}$. It turns out that $Q$ is the inverse of
$i_*$ after localizing at $(f)$:  $Q \circ i_* =1$ because of
lemma \ref{pullbackpushforward}, and $i_* \circ Q =1$ because the
module structure is compatible with projections; namely for $\pi :
\nu \to F$ and $a \in H^\bullet_{G}(M;\alpha)$, one has that
$\pi^* (i^* a) \wedge \pi^* (e_G(\nu)^{-1}) = \pi^*(i^*a \wedge
e_G(\nu)^{-1})$. So, we can conclude:

\begin{theorem} [Localization at fixed points] For all $x \in
H^\bullet_{G}(M;\alpha)$ in a suitable localization, one has
$$x = \sum_{Z \subset F} i_*^Z (i_Z^*( x)) \wedge e_G(\nu_Z)^{-1}$$
\end{theorem}

From the localization theorem we get that $H^\bullet_G(\TTC M)
_{(f)} \cong H^\bullet_G(\TTC F)_{(f)}$, and from the results in
section \S\ref{equicoho:trivial} we get that $H^\bullet_G(\TTC F)$
is isomorphic to the twisted equivariant cohomology of the fixed
point set. So, for $\alpha = H + \sum_j u_j \xi_j$ the twisting
form in $M$, we have
\begin{corollary}
If the $G$ action is pure, then $$H^\bullet_G(\TTC M)_{(f)} \cong H^\bullet_{G}(F;i^*\alpha)_{(f)}$$ and if the 1-forms $i^*\xi_j$ are all exact (say when $H^1(F)=0$) then
$$H^\bullet_G(\TTC M)_{(f)} \cong H^\bullet(F, i^*H) \tensor \widehat{H}^\bullet(BG)_{(f)}=  H^\bullet(F, i^*H) \tensor
\CC[[u_1, \dots, u_n]]_{(f)}.$$
\end{corollary}

\begin{example} {\rm
Let's consider the $S^1$ extended action on $S^3$ from section \S
\ref{section:non-free_action_example}. A point in $S^3$ is a pair
of complex numbers $z=(z_1,z_2)$ with $|z_1|^2 + |z_2|^2 =1$ that
could also be written in  polar coordinates as $z_1 = e^{i \phi_1}
\sin \lambda$ and $z_2 = e^{i \phi_2} \cos \lambda$. The 3-form $H
$ is $-\sin (2\lambda) d \lambda \wedge d \phi_1 \wedge d \phi_2
$,  the $S^1$ action is defined by rotating the first coordinate
$z_1$ and the extended action is $\RR^1 \to \Gamma(\TT S^3) : 1
\mapsto \frac{\partial}{\partial \phi_1} - \cos^2\lambda d\phi_2$.
As the action is pure, the form $\alpha := H - u \cos^2\lambda
d\phi_2$ is equivariantly closed, and therefore defines a
cohomology class $[\alpha] \in H^3_{S^1}(S^3)$. Recall that
$H^\bullet(BS^1) = \CC[u]$.

The fixed point set of the circle action is the set $F =
\{(z_1,z_2) | z_1=0\} \cap S^3$ which is also a circle.
If $i : F \to S^3$ is the inclusion, then $i^* \alpha =-u d\phi_2 \in H^3_{S^1}(F) = H^1(F) \otimes H^2(BS^1)$ and therefore we can apply the results of the example \ref{example:trivial_action}. So we have that $H^{od}_{S^1}(F;i^* \alpha) =  \RR$, while $H^{ev}_{ S^1}(F;i^* \alpha) = 0$. 
As the normal bundle $\nu$ of $F$ in $S^3$ is trivial, and the action of the circle in the fibers is by rotation, then the equivariant euler class of the normal bundle is $e_{S^1}(\nu) = u$. By the localization theorem, if we invert  $u$ we get the isomorphisms
$$H^{od}_{ S^1}(S^3;\alpha)_{(u)} \cong H^{od}_{ S^1}(F;i^*\alpha)_{(u)} \cong
\RR_{(u)}=0, $$ 
$$\mbox{and} \ \ \ H^{ev}_{S^1}(S^3;\alpha)_{(u)} \cong H^{ev}_{ S^1}(F;i^*\alpha)_{(u)}=0.$$
Hence, applying the localization theorem we can deduce that the equivariant twisted cohomology $H^\bullet_{ S^1}(S^3;\alpha)$ is a torsion $\CC [u]$-module.

}
\end{example}

\bibliographystyle{amsplain}
\bibliography{equiv_extended}

\end{document}